\documentclass[11pt]{article}
\usepackage{amsfonts,amsmath,amsxtra}
\usepackage{amssymb,amscd}
\usepackage[matrix,arrow]{xy}
\usepackage{latexsym}

\usepackage[T2A]{fontenc}
\usepackage[english,russian]{babel}
\usepackage[utf8]{inputenc}

\def\hybrid{\topmargin 0pt      \oddsidemargin 0pt
        \headheight 0pt \headsep 0pt
        \textwidth 16.5cm
        \textheight 23cm
        \marginparwidth 0.0in
        \parskip 5pt plus 1pt   \jot = 1.5ex}
\catcode`\@=11
\def\marginnote#1{}
\newcount\hour
\newcount\minute
\newtoks\amorpm
\hour=\time\divide\hour by60 \minute=\time{\multiply\hour by60
\global\advance\minute by-\hour}
\edef\standardtime{{\ifnum\hour<12 \global\amorpm={am}%
        \else\global\amorpm={pm}\advance\hour by-12 \fi
        \ifnum\hour=0 \hour=12 \fi
      \number\hour:\ifnum\minute<10 0\fi\number\minute\the\amorpm}}
\edef\militarytime{\number\hour:\ifnum\minute<10 0\fi\number\minute}

\def\draftlabel#1{{\@bsphack\if@filesw {\let\thepage\relax
   \xdef\@gtempa{\write\@auxout{\string
      \newlabel{#1}{{\@currentlabel}{\thepage}}}}}\@gtempa
   \if@nobreak \ifvmode\nobreak\fi\fi\fi\@esphack}
        \gdef\@eqnlabel{#1}}
\def\@eqnlabel{}
\def\@vacuum{}
\def\draftmarginnote#1{\marginpar{\raggedright\scriptsize\tt#1}}

\def\draft{\oddsidemargin -0.1truein
        \def\@oddfoot{\sl preliminary draft \hfil
        \rm\thepage\hfil\sl\today\quad\militarytime}
        \let\@evenfoot\@oddfoot \overfullrule 3pt
        \let\label=\draftlabel
        \let\marginnote=\draftmarginnote
\def\@eqnnum{{\rm (\theequation)}
\rlap{\kern\marginparsep\tt\@eqnlabel}%
\global\let\@eqnlabel\@vacuum}  }
\newfont{\Bbbb}{msbm7 scaled 1\@ptsize00}
\newcommand{\zs}{\raise-1pt\hbox{$\mbox{\Bbbb Z}$}}

scaled 1\@ptsize00
\font\sevenmsa=msam6 %scaled 1\@ptsize00
\newfam\msafam
\textfont\msafam=\sevenmsa
\def\hexnumber@#1{\ifnum#1<10 \number#1\else
\ifnum#1=10 A\else\ifnum#1=11 B\else\ifnum#1=12 C\else \ifnum#1=13
D\else\ifnum#1=14 E\else\ifnum#1=15 F\fi\fi\fi\fi\fi\fi\fi}
\def\msa@{\hexnumber@\msafam}
\def\llcorner{\delimiter"4\msa@78\msa@78 }
\def\lrcorner{\delimiter"5\msa@79\msa@79 }
\mathchardef\blacktriangleright="3\msa@49
\mathchardef\blacktriangleleft="3\msa@4A \font\tenmsb=msbm10 scaled
1\@ptsize00
\newfam\msbfam
\textfont\msbfam=\tenmsb \scriptfont\msbfam=\tenmsb
\newdimen\Squaresize \Squaresize=14pt
\newdimen\Thickness \Thickness=0.5pt

\def\Square#1{\hbox{\vrule width \Thickness
   \vbox to \Squaresize{\hrule height \Thickness\vss
      \hbox to \Squaresize{\hss#1\hss}
   \vss\hrule height\Thickness}
\unskip\vrule width \Thickness} \kern-\Thickness}

\def\Vsquare#1{\vbox{\Square{$#1$}}\kern-\Thickness}

\def\numberbysection{\@addtoreset{equation}{section}
        \def\theequation{\thesection.\arabic{equation}}}
\numberbysection

\renewcommand{\theequation}{\thesection.\arabic{equation}}
\def\titlepage{\@restonecolfalse\if@twocolumn\@restonecoltrue\onecolumn
     \else \newpage \fi \thispagestyle{empty}\c@page\z@
        \def\thefootnote{\fnsymbol{footnote}} }

\def\endtitlepage{\if@restonecol\twocolumn \else  \fi
        \def\thefootnote{\arabic{footnote}}
        \setcounter{footnote}{0}}  %\c@footnote\z@ }
\relax

\hybrid
\makeatletter
\newdimen\normalarrayskip            % skip between lines
\newdimen\minarrayskip               % minimal skip between lines
\normalarrayskip\baselineskip \minarrayskip\jot
\newif\ifold             \oldtrue            \def\new{\oldfalse}
\def\arraymode{\ifold\relax\else\displaystyle\fi}%mode of array enrties
\def\eqnumphantom{\phantom{(\theequation)}} % ight phantom in eqnarray
\def\@arrayskip{\ifold\baselineskip\z@\lineskip\z@
     \else
     \baselineskip\minarrayskip\lineskip1\baselineskip\fi}

\def\@arrayclassz{\ifcase \@lastchclass \@acolampacol \or
\@ampacol \or \or \or \@addamp \or
   \@acolampacol \or \@firstampfalse \@acol \fi
\edef\@preamble{\@preamble
  \ifcase \@chnum
     \hfil$\relax\arraymode\@sharp$\hfil
     \or $\relax\arraymode\@sharp$\hfil
     \or \hfil$\relax\arraymode\@sharp$\fi}}

\def\@array[#1]#2{\setbox\@arstrutbox=\hbox{\vrule
     height\arraystretch \ht\strutbox
     depth\arraystretch \dp\strutbox
width\z@}\@mkpream{#2}\edef\@preamble{\halign \noexpand\@halignto
\bgroup \tabskip\z@ \@arstrut \@preamble \tabskip\z@ \cr}%
\let\@startpbox\@@startpbox \let\@endpbox\@@endpbox
  \if #1t\vtop \else \if#1b\vbox \else \vcenter \fi\fi
  \bgroup \let\par\relax
  \let\@sharp##\let\protect\relax
  \@arrayskip\@preamble}
\def\eqnarray{\stepcounter{equation}%
              \let\@currentlabel=\theequation
              \global\@eqnswtrue
              \global\@eqcnt\z@
              \tabskip\@centering              %formulae  centering
              \let\\=\@eqncr
              $$%
            \halign to \displaywidth  \bgroup
             \eqnumphantom \@eqnsel
      \hskip\@centering                               %right tab%
    $\displaystyle  \tabskip\z@ {##}$%
    &\global\@eqcnt\@ne \hskip 2\arraycolsep
         $ \displaystyle  \arraymode{##}$\hfil
    &\global\@eqcnt\tw@ \hskip 2\arraycolsep
         $\displaystyle\tabskip\z@{##}$\hfil
         \tabskip\@centering
    &{##}\tabskip\z@\cr}
\makeatother
\def\CH {\mathcal{H}}

\def\CQ {\mathcal{Q}}

\def\pr {\partial}

\newtheorem{prop}{Proposition}[section]           %  ETC ...

\newtheorem{lem}{Lemma}[section]

\newtheorem{conj}{Conjecture}[section]

\newcommand\bqa{\begin{eqnarray}}
\newcommand\eqa{\end{eqnarray}}
\def\be{\begin{eqnarray}\new\begin{array}{cc}}
\def\ee{\end{array}\end{eqnarray}}

\def\beq{\begin{equation}}
\def\eeq{\end{equation}}
\def\bse{\begin{subequations}}                %%%SUBEQUATIONS
\def\ese{\end{subequations}}
\def\bp{\begin{pmatrix}}
\def\ep{\end{pmatrix}}

\def\ux{\underline{x}}
\def\uy{\underline{y}}
\def\uz{\underline{z}}

\def\stack#1#2{\raise0.7pt\hbox{$\mathrel{\mathop{#2}\limits^{#1}}$}}
\def\tr{\triangleright}
\def\tl{\triangleleft}
\def\sem{\mathsurround=0pt \raise1pt
\hbox{$\scriptscriptstyle>\!\!$}\:\!\!\tl}
\def\mes{\mathsurround=0pt \tr\!\:\!\raise0.8pt
\hbox{$\scriptscriptstyle\!\!<$}\,}
\def\]{\mathsurround=0pt ]\raise-2pt\hbox{$_\ast$}}
\def\<{\langle}
\def\>{\rangle}

\def\CQ{{\cal Q}}

\def\CH{\mathcal{H}}

\def\we{\raise-1pt\hbox{$\,\stackrel{\wedge}{,}\,$}}
\def\tr{{\rm tr}\,}

\def\pr {\partial}

\newcounter{pac}[section]
\newcounter{pacc}[subsection]
\newcommand{\npa}{\addtocounter{pacc}{1} \noindent {\bf
\arabic{section}.\arabic{subsection}.\arabic{pacc}}\,\,\,}
 0

 0

\begin{document}

%\title{\bf Quantum Toda Chains Intertwined}
%\author{Anton Gerasimov, Dimitri Lebedev, and Sergey Oblezin}
%\date{}
%\maketitle

\begin{flushright}
\selectlanguage{russian}
``Наиболее удобной моделю для отработки \\
 такой связи
является цепочка Тоды'' \\
\vspace{3 mm}
Л.Д. Фаддеев, Препринт ЛОМИ Р-2-79
\end{flushright}
\selectlanguage{english}

\vspace{15 mm}
\centerline{\Large \bf Quantum Toda Chains Intertwined}
\vspace{4 mm}
\centerline{Anton Gerasimov, Dimitri Lebedev, Sergey Oblezin}
\vspace{7 mm}
%\vspace{-3 mm}

\centerline{\it To Ludvig Dmitrievich Faddeev on the occasion of his
  75-th birthday}
\vspace{5 mm}
\begin{abstract}
\noindent {\bf Abstract}. We conjecture an explicit construction of
integral operators intertwining various quantum Toda chains.
Compositions of the intertwining operators provide recursive and
$\CQ$-operators for quantum Toda chains.
In particular we propose a generalization
of our previous  results for Toda chains corresponding
to classical Lie algebra  to the generic $BC_n$ and Inozemtsev Toda chains.
We also conjecture explicit form of $\CQ$-operators for closed Toda chains
corresponding to Lie algebras $B_{\infty}$, $C_{\infty}$,
$D_{\infty}$, affine Lie algebras $B^{(1)}_n$, $C^{(1)}_n$,
$D^{(1)}_n$, $D^{(2)}_n$, $A^{(2)}_{2n-1}$,  $A^{(2)}_{2n}$ and
the affine analogs of $BC_n$ and Inozemtsev Toda chains.
\end{abstract}

\vspace{1cm}

%\end{titlepage}

%\draft                             %SWITCH ON/OFF DRAFT VERSION%

%\tableofcontents
\normalsize
\section*{Introduction}

An interesting  integral representation  for common
eigenfunctions of the  quantum $A_n$-Toda chain  Hamiltonian operators was
proposed in  \cite{Gi}. Using the representation theory approach to Toda chain
eigenfunctions the representation \cite{Gi} was generalized to other
classical series $B_n$, $C_n$ and $D_n$ \cite{GLO2}. In \cite{GKLO}, \cite{GLO2}
we stress an important role of the elementary intertwining
operators in the construction of solutions of
quantum Toda chains.  The common eigenfunctions of quantum
Toda Hamiltonians corresponding to the same classical series
can be obtained  by a recursive  application of some integral
operators. In general the recursive integral operators can be  naturally
represented as a product of elementary integral operators
intertwining quantum Toda Hamiltonians corresponding to different
classical series.
The integral kernels of the elementary
intertwining  operators have a uniform structure of the  exponent of
a linear combination of exponents in natural coordinates.
This decomposition of the recursive operators into elementary
intertwiners seems has a methodological importance
for quantization of classical Toda chains.
 Indeed, for classical Toda chains  the  recursive operators
 are given by complicated canonical transformations. The direct quantization of
these canonical transformations  is far from obvious.
Let classical canonical
transformations be naturally split into  products of elementary canonical
transformations allowing standard quantization. Then
the quantization of the canonical transformation is easily obtained
 as a product of quantized elementary canonical transformations.
As it was argued  in \cite{GLO2} and is demonstrated in this paper
this is precisely the case for Toda chains.

One can expect that the algebra of intertwining operators acting
between various Toda chains should have an interesting
interpretation in representation theory of classical series of Lie
groups. The representation theory interpretation of a particular
instance of the integral intertwining operator commuting with
quantum Toda chain Hamiltonians  was  proposed in
\cite{GLO3}. These operators were  identified with the generators of
Hecke algebras of the classical Lie groups with respect to maximal
compact subgroups. The affine version of these operators were
introduced long ago by Pasquier and Gaudin as $\CQ$-operator for
affine Toda chain \cite{PG}.

In this paper  we extend the integral intertwiner
operator approach to open Toda chains whose representation theory
interpretation for generic coupling constants is not known.
These include $BC_n$-Toda chains 
constructed by Sklyanin  \cite{S} and Toda  chains
proposed by  Inozemtsev in the classical setting  \cite{I}.
In the following we refer to the latter series of Toda chains
as $I_n$-Toda chains.
We propose elementary intertwining operators for $BC_n$ and $I_n$
Toda chains and check the intertwining properties with respect to
quadratic quantum Hamiltonians.  The elementary intertwiners appear to
relate $BC$ and $I$-Toda chains thus stressing  the importance of the Inozemtsev
Toda chain series as an honorable member of the Toda chains family.
The constructed elementary intertwiners allow us to propose
explicit integral expressions for common eigenfunctions of
$BC_n$- and $I_n$-Toda chain Hamiltonians for generic spectrum.

In the second part of the paper we propose analogs of the
Pasquier-Gaudin integral $\CQ$-operator \cite{PG} for affine Toda
chains based on classical series of affine Lie algebras $A_n^{(1)}$,
$A^{(2)}_{2n}$, $A_{2n-1}^{(2)}$, $B_n^{(1)}$, $C_n^{(1)}$,
$D_n^{(1)}$, $D_n^{(2)}$  infinite Lie algebras $B_{\infty}$,
$C_{\infty}$ and $D_{\infty}$, affine and semi-infinite version of
$BC_n$ and $I_n$ Toda-chains.

Finally let us stress that the explicit integral expressions for
common eigenfunctions constructed using elementary intertwiners
might provide a general approach to the  quantization of classical
integrable systems. This program was completely realized for
$A_n$-Toda chains in \cite{GLO3}. The first step is a
construction of recursive operators as a product of quantized
elementary intertwining integral operators and checking commutation
relations with quadratic quantum Hamiltonian. The next step is a
construction of quadratic Hamiltonian eigenfunctions using the
recursive operators. Finally one should check orthogonality and
completeness relations for the constructed eigenfunctions.
Note that the complete set of eigenfunctions can be considered
(after a standard Fourier transform)
as an integral kernel of a unitary transformation $U$ intertwining
quantum Toda chain with the quantum theory of free particles
(i.e. having zero potential).  The higher quantum Hamiltonians can now
be obtained by a conjugation by $U$ of  the
standard set of the higher quantum Hamiltonians for a collection of
free particles. The correct classical limit of these Hamiltonians is
guaranteed by the correct classical limit of the quantum recursive
operators used to construct $U$. It would be interesting to apply this
line of reasoning to  $BC_n$- and $I_n$-Toda chains
which do not allow simple representation theory realization.
What is missed at the moment are the orthogonality and completeness
relations for these sets of  functions.  The proof
of the orthogonality and completeness relations consists of routine
manipulations with iterative integrals and will be published
elsewhere.  Modern approach to solving quantum integrable systems
based on the Quantum Inverse Scattering Method (QISM)   achieved
a spectacular success for various classes of quantum integrable
models (see \cite{F} for a concise account). One might hope that
the approach based on the explicit integral operators used in this
paper might add some complimentary information about the structure of
quantum integrable systems.

This paper is an  updated version of the earlier note \cite{GLO1}.
After \cite{GLO1}, \cite{GLO2} the paper \cite{KS} with the
results on classical recursion operators for $BC_n$-Toda chains
appeared.  We expect that the classical limit of our constructions
for $BC_n$ Toda chain  should be compatible with the results of \cite{KS}.

The plan of this paper  is as follows. In Section 1 we summarize the
results of \cite{Gi}, \cite{GKLO}, \cite{GLO2} for the case of zero
spectrum. In Section 2 we propose elementary integral intertwining
operators and integral representations for $BC_n$- and $I_n$-Toda
chain eigenfunctions. In Section 3 analogs of Pasquier-Gaudin
$\CQ$-operator at zero spectral parameter are proposed for affine
versions of Toda chains corresponding to classical Lie groups and
$BC_n$- and $I_n$-Toda chains.
 In Section 4 the analogous conjectures are stated for semi-infinite Toda
 chains.

{\em Acknowledgments}: The research was supported by grants
RFBR-08-01-00931-a, \selectlanguage{russian}
09-01-93108-НЦНИЛ-а\selectlanguage{english}, 
 AG was  also partly supported by
Science Foundation Ireland grant.

\section{Integral representation for
the eigenfunctions of quantum Toda chains for classical groups}

We start reviewing  a construction of the intertwining
 operators and common eigenfunction of
quantum Toda chains corresponding to classical Lie
algebras \cite{GLO2}. Then  we conjecture  explicit
integral representations for eigenfunctions for
generic $BC_n$ and $I_n$ Toda chains.

Let us first recall the standard facts on  Toda chains
 corresponding to arbitrary root systems (see e.g. \cite{RSTS}).
Let $\mathfrak{g}$ be a simple Lie algebra, $\mathfrak{h}$ be a
Cartan subalgebra, $n=\dim \mathfrak{h}$ be the rank of
$\mathfrak{g}$, $R\subset \mathfrak{h}^*$  be the root system,  $W$
be  the  Weyl group. Let us fix a decomposition $R=R_+\cup R_-$ of
the roots on positive and negative roots. Let
$\alpha_1,\ldots,\alpha_n$ be the bases of  simple roots. Let $(,)$
be a $W$-invariant bilinear symmetric form on $\mathfrak{h}^*$
normalized so that $(\alpha,\alpha)=2$ for a long root. This form
provides an identification of $\mathfrak{h}$ with $\mathfrak{h}^*$
and thus can be considered as a bilinear form on $\mathfrak{h}$.
Choose an orthonormal basis $e=\{e_1,\ldots,e_{n}\}$ in
$\mathfrak{h}$. Then for any $x\in \mathfrak{h}$ one has  a
decomposition $x=\sum_{i=1}^n x_ie_i$. One associates  with these
data an open Toda chain with a quadratic Hamiltonian
 \bqa\label{Hamtwo}
  H_2^{R}(\ux)= -\frac{1}{2}\sum_{i=1}^n
  \frac{\partial^2}{\partial x_i^2}\,+\,
  \sum_{i=1}^ng_ie^{\alpha_i(x)}\,
 \eqa
where $g_i$ are (in general complex) coupling constants.  
This Hamiltonian is a one of the generators $\CH_i^R$,
$i=1,2,\ldots, {\rm rank}(R)$ of the commutative ring of quantum
Hamiltonians of the corresponding Toda chain. The common
eigenfunctions of the quantum Hamiltonians satisfy the equations
$$
H_i^{R}(\ux)
\Psi_{\underline{\lambda}}^{R}(\ux)=
c_i(\underline{\lambda})\,
 \Psi^{R}_{\underline{\lambda}}(\ux),
$$
where $c_i(\underline{\lambda})$ are  polynomials of
the complex variables $\underline{\lambda}=(\lambda_1,\ldots ,\lambda_n)$
parameterizing spectrum.  In the following we
will consider only the eigenfunctions  corresponding to zero
eigenvalues of the Hamiltonians and omit the index
$\underline{\lambda}$. Borrowing Physics terminology we will
call such eigenfunctions ground state wave
functions.

\subsection{Integral representation of $A_{n-1}$-Toda chain eigenfunctions}

In this Subsection  we recall a recursive construction \cite{GKLO} of the Givental
integral representation of $\mathfrak{gl}_{n}$-Toda chain
eigenfunctions \cite{Gi}. Quantum
$\mathfrak{gl}_{n}$-Toda chain is characterized  by its
linear and quadratic Hamiltonians
 \be\label{Aham}
  H_1^{\mathfrak{gl}_n}(\ux)\,=\,\,\imath \,\,
\sum\limits_{i=1}^n
\, \frac{\partial}{\partial x_i}\,,\\
  H_2^{\mathfrak{gl}_n}(\ux)\,=\,-\frac{1}{2}\sum\limits_{i=1}^n
  \frac{\partial^2}{\partial x_i^2}\,+\,
  \sum\limits_{i=1}^{n-1}g_i e^{x_{i+1}-x_{i}}\,,
 \ee
where $\ux=(x_1,\ldots,x_n)$ and $g_i$ are coupling constants. These
Hamiltonians are the first two of a series of  generators $H_i^{\mathfrak{gl}_{n}}$,
$i=1,2,\ldots n$  of the commutative ring of quantum Hamiltonians of
$\mathfrak{gl}_{n}$-Toda chain. The common eigenfunctions of the
quantum Hamiltonians satisfy the equations
\be\label{eigencond}
H_i^{\mathfrak{gl}_n}(\ux)
\Psi_{\underline{\lambda}}^{\mathfrak{gl}_n}(\ux)=
c_i(\underline{\lambda})\,
 \Psi^{\mathfrak{gl}_n}_{\underline{\lambda}}(\ux),
\ee
where $c_i(\underline{\lambda})$ are the symmetric polynomials
of the complex variables $\underline{\lambda}=(\lambda_1,\ldots
,\lambda_n)$ parameterizing spectrum.  For instance
$c_1(\underline{\lambda})=\lambda_1+\ldots +\lambda_n$ and
$c_2(\underline{\lambda})=\lambda_1^2+\ldots+
\lambda_n^2$. Common eigenfunctions for $\mathfrak{sl}_n$-Toda chain
are obtained by imposing the additional condition
$c_1(\underline{\lambda})=0$ in \eqref{eigencond}. In the following
we will consider ground state wave functions
corresponding to $\underline{\lambda}=0$
and thus will not  distinguish $\mathfrak{gl}_n$- and
$\mathfrak{sl}_n$-Toda chains.

\npa The following integral representation
 for a  ground state wave function  quantum Hamiltonians
of $\mathfrak{gl}_n$-Toda chain was derived in \cite{Gi}
 \be\label{intrep}
\Psi^{\mathfrak{gl}_n}(\ux)=\int_{\Gamma}
e^{\mathcal{F}_n(x)}\prod_{k=1}^{n-1}\prod_{i=1}^k
dx_{k,i}, \ee
where $x_{n,i}:=x_i$, $\ux=(x_1,\ldots,x_{n})$,  the function $\mathcal{F}_n(x)$
is given by
 \be\label{pot}
  \mathcal{F}_n(x)\,=\,-\sum_{k=1}^{n-1}\sum_{i=1}^k
  \Big(e^{x_{k,i}-x_{k+1,i}}+g_ie^{x_{k+1,i+1}-x_{k,i}}\Big),
 \ee
and the integration domain $\Gamma$ is a middle dimensional
submanifold in the $n(n-1)/2$- dimensional complex torus with
coordinates $\{\exp\,{x_{k,i}},\, i\leq 1\leq k\leq n-1\}$. The
choice of the integration contour is dictated by the condition
on the integrand to decrease exponentially   at infinities.

The integral representation (\ref{intrep}) can be derived using the
recursion over the rank  $n$ of the Lie algebra $\mathfrak{gl}_n$ \cite{GKLO}.
Indeed (\ref{intrep}) is naturally represented in the recursive form
 \bqa\label{recintrep}
  \Psi^{\mathfrak{gl}_n}(\ux)\,=\,
  \int_{\Gamma}\prod_{k=1}^{n-1}
  Q_{k+1,\,k}(\ux_{k+1};
  \,\ux_k)\,d\ux_k\,,
 \eqa with the integral kernel
 \be\label{QBAXTER}
  Q_{k+1,k}(\ux_{k+1};\,\ux_k)\,=\,\exp\Big\{\,-\sum_{i=1}^{k}
  \Big(e^{x_{k,i}-x_{k+1,i}}+g_ie^{x_{k+1,i+1}-x_{k,i}}\Big)\,\Big\}.
 \ee
Here we define $\ux_k=(x_{k,1}\,,\ldots,x_{kk})$ and
$x_i:=x_{n,i},\,1\leq i\leq n$. The zero eigenvalue
property  of (\ref{intrep}) with respect to quantum Hamiltonians
is then translated into the intertwining
relations between integral operators with the kernels   $Q_{k+1,k}$
and Hamiltonians $H_i^{\mathfrak{gl}_k}$ and
$H_i^{\mathfrak{gl}_{k+1}}$. In particular for the quadratic
Hamiltonian \eqref{Aham} we have the following easily verified
property
 \bqa\label{intertw}
  H^{\mathfrak{gl}_{k+1}}(\ux_{k+1})Q_{k+1,k}(\ux_{k+1},\,
  \ux_k)\,=\,
  Q_{k+1,k}(\ux_{k+1},\, \ux_k)\,H^{\mathfrak{gl}_k}(\ux_k),
 \eqa where
 \be
  H^{\mathfrak{gl}_k}(\ux_k)\,=\,-\frac{1}{2}\sum\limits_{i=1}^k
  \frac{\partial^2}{\partial x_{k,i}^2}\,+\,
  \sum\limits_{i=1}^{k-1}g_i e^{x_{k,i+1}-x_{k,i}}\, .
 \ee
Here and in the following we assume that in the relations similar to
(\ref{intertw}) the Hamiltonian operator on l.h.s. acts on the right
and the Hamiltonian on r.h.s. acts on the left. From \eqref{intertw}
one can easily deduce
\bqa
\Big(\,-\frac{1}{2}\sum\limits_{i=1}^n
  \frac{\partial^2}{\partial x_i^2}\,+\,
  \sum\limits_{i=1}^{n-1}g_i e^{x_{i+1}-x_{i}}\,\Big)\,\Psi^{\mathfrak{gl}_n}(\ux)= 0.
 \eqa
The eigenvalue property for higher Hamiltonians is not so easy to
verify and deserve more  involved considerations
\cite{Gi}, \cite{GKLO}.

\npa The integral operator defined by the kernel (\ref{QBAXTER}) is
closely related with a $\CQ$-operator realizing B\"{a}cklund
transformations in a  Toda chain corresponding to affine Lie algebra
$\widehat{\mathfrak{gl}}_n$. This is a particular example of
the integral intertwining operator introduced by Pasquier and Gaudin
\cite{PG}.  The integral kernel is given by
 \bqa\label{baff}
  Q^{\widehat{\mathfrak{gl}}_n}(\ux,\uy)=\exp\Big\{\,-
  \sum_{i=1}^n\Big(e^{x_{i}-y_{i}}+g_ie^{y_{i+1}-x_i}\Big)\,\Big\},
 \eqa
where $x_{i+n}=x_i$ and $y_{i+n}=y_i$. This operator commutes with
the Hamiltonian operators of the closed Toda chain. Thus for example
for the quadratic Hamiltonian we have
 \bqa\label{intertwaff}
  H_2^{\widehat{\mathfrak{gl}}_n}(\underline{x})Q^{\widehat{\mathfrak{gl}}_n}(\ux,\,
  \uy)\,=\,Q^{\widehat{\mathfrak{gl}}_n}(\ux,\,\uy)
  H_2^{\widehat{\mathfrak{gl}}_n}(\uy),
 \eqa where
 \be
  H_2^{\widehat{\mathfrak{gl}}_n}(\ux)\,=\,
  -\frac{1}{2}\sum\limits_{i=1}^{n}\frac{\partial^2}{\partial
  x_i^2}\,+\,
  \sum\limits_{i=1}^{n}g_i e^{x_{i+1}-x_{i}}\,.
 \ee Here we impose the
conditions  $x_{i+n}=x_i$. The recursive operator (\ref{QBAXTER})
can be obtained from the operator (\ref{baff}) in the limit
$g_{n}\rightarrow 0$, $x_n\rightarrow-\infty$ in such a way that
$g_ne^{-x_n}\to 0$.

\subsection{Toda chains corresponding to classical Lie algebras}

 In this Subsection we recall the explicit form of the quadratic
quantum Hamiltonians of Toda chains corresponding to classical Lie
algebras.

\npa Let $e=\{e_1,\ldots,e_{n}\}$ be an  orthonormal basis in
$\mathbb{R}^n$. The  root system of type $B_n$ can be defined as
 \bqa
  \alpha_1=e_1,\hspace{1.5cm}\alpha_{i+1}=e_{i+1}-e_i,\hspace{2cm}
  1\leq i\leq n-1,
 \eqa
and the corresponding Dynkin diagram is
 \bqa
 \xymatrix{
  \alpha_1\ar@{<=}[r] & \alpha_2 & \ldots\ar@{-}[l] &
  \alpha_{n-1}\ar@{-}[l]\ar@{-}[r] & \alpha_n
 }\,.
 \eqa
Quadratic Hamiltonian operator of the corresponding Toda chain
is given by
 \bqa
  H_2^{B_n}(\ux)\,=\,-\frac{1}{2}\sum_{i=1}^{n}
  \frac{\partial^2}{\partial x_i^2}\,+\,
  g_1e^{x_1}\,+\,\sum_{i=1}^{n-1}g_{i+1}e^{x_{i+1}-x_i}\,.
 \eqa

\npa The  root system of type $C_n$ can be defined as
 \bqa
  \alpha_1=2e_1,\hspace{1.5cm} \alpha_{i+1}=e_{i+1}-e_i,\hspace{2cm}
  1\leq i\leq n-1,
 \eqa
and the   corresponding Dynkin diagram is
 \bqa
 \xymatrix{
  \alpha_1\ar@{=>}[r] & \alpha_2 & \ldots\ar@{-}[l] &
  \alpha_{n-1}\ar@{-}[l]\ar@{-}[r] & \alpha_n
 }\,.
 \eqa
Quadratic Hamiltonian operator of the corresponding Toda chain
is given by
 \be
  H^{C_n}(\uz)\,=\,-\frac{1}{2}\sum_{i=1}^n
  \frac{\partial^2}{\partial z_i^2}\,+\,2g_1e^{2z_1}\,+\,
  \sum_{i=1}^{n-1}g_{i+1}e^{z_{i+1}-z_i}\,, \ee
where $\uz=(x_1,\ldots,z_n)$.

There is a more general Toda chain associated with
root system $C_n$ (see e.g. \cite{RSTS}, Remark p.61)
with the quadratic Hamiltonian
 \bqa\label{BCstar}
  H_2^{BC^*_n}(\uz)\,=\,-\frac{1}{2}\sum_{i=1}^n
  \frac{\partial^2}{\partial z_i^2}\,-\,
  \frac{g_1}{2}e^{z_1}+\frac{g_1^2}{2}e^{2z_1}\,+\,
  \sum_{i=1}^{n-1}g_{i+1}e^{z_{i+1}-z_i}\,,
 \eqa

The $BC_n^*$-Toda chain can be obtained
as a degeneration of the generic $BC_n$-Toda chain discussed in the
next. This explains our notations.

\npa The root system of type $D_n$ is
 \bqa
  \alpha_1=e_1+e_2,\hspace{1.5cm}
  \alpha_{i+1}=e_{i+1}-e_i,\hspace{2cm}1\leq i\leq n-1,
 \eqa and the
corresponding Dynkin diagram is
$$
 
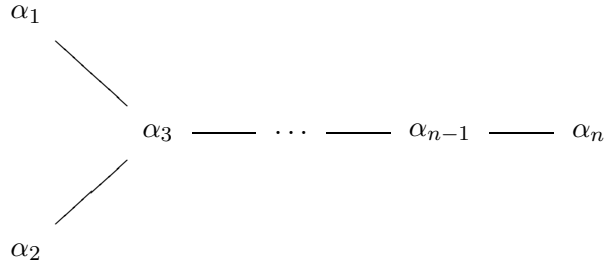
\xymatrix{
 \alpha_1\ar@{-}[dr] &&&&&\\
 & \alpha_3\ar@{-}[r] & \ldots\ar@{-}[r] &
 \alpha_{n-1}\ar@{-}[r] & \alpha_n\\
 \alpha_2\ar@{-}[ur] &&&&&
}\,\,\,
$$
Quadratic Hamiltonian operator of the $D_n$-Toda chain is given by
 \be
  H^{D_n}(\ux)\,=\,-\frac{1}{2}\sum_{i=1}^n
  \frac{\partial^2}{\partial x_i^2}\,+\,
  g_1g_2e^{x_1+x_2}\,+\,\sum_{i=1}^{n-1}g_{i+1}e^{x_{i+1}-x_i}\,.
 \ee

\subsection{Intertwining integral operators}

 We start with the elementary intertwining operators relating Toda
chains corresponding to different series. These elementary
intertwiners have simple integral kernels given by exponents
of the linear combination of exponents. The integration contours 
 are chosen in such way that the integrals considered below converge. 

\npa Integral kernel
 \be\label{QBC1}
  Q^{BC^*_n}_{\,\,\,\,B_{n-1}}(\uz_n;\,\ux_{n-1})\,=\,
  \exp\Big\{-\Big(\,g_1e^{z_{n,1}}+\sum_{i=1}^{n-1}\Big(
  e^{x_{n-1,i}-z_{n,i}}+g_{i+1}e^{z_{n,i+1}-x_{n-1,i}}\Big)\,\Big)\Big\},
 \ee
satisfies the following relations
 \bqa\label{BCB}
  H^{BC^*_n}(\uz_n)\,Q^{BC^*_n}_{\,\,\,\,B_{n-1}}(\uz_n,\,\ux_{n-1})\,=\,
  Q^{BC^*_n}_{\,\,\,\,B_{n-1}}(\uz_n,\,\ux_{n-1}) \,H^{B_{n-1}}(\ux_{n-1})\,.
 \eqa
Thus the corresponding integral operator
intertwines quadratic Hamiltonians of $BC^*_n$- and $B_{n-1}$-Toda chains.
Similarly an elementary operator intertwining quadratic $B_n$ and $BC^*_{n}$
 Hamiltonians  has an integral kernel
 \be\label{QBC2}
  Q^{B_n}_{\,\,\,\,BC^*_n}(\ux_n;\,\uz_n)\,=\\
  \exp\Big\{-\Big(\,g_1e^{z_{n,1}}+\sum_{i=1}^{n-1}\Big(e^{x_{n,i}-z_{n,i}}+
  g_{i+1}e^{z_{n,i+1}-x_{n,i}}\Big)+e^{x_{nn}-z_{nn}}\,\Big)\Big\}.
 \ee

\npa An integral operator intertwining $C_n$ and $D_{n}$ Hamiltonians has
the kernel
 \be\label{QC1}
 \hspace{-5mm}
  Q^{C_n}_{\,\,\,\,D_n}(\uz_n;\,\ux_n)\,=\,
  \exp\Big\{-\Big(\,g_1e^{x_{n,1}+z_{n,1}}+
  \sum_{i=1}^{n-1}\Big(e^{z_{n,i}-x_{n,i}}+
  g_{i+1}e^{x_{n,i+1}-z_{n,i}}\Big)+e^{z_{nn}-x_{nn}}\Big)\Big\}.
\ee
Similarly an integral operator with the kernel
 \be\label{QC2}
  Q^{D_n}_{\,\,\,C_{n-1}}(\ux_n;\,\uz_{n-1})\,=\,
  \exp\Big\{-\Big(\,g_1e^{x_{n,1}+z_{n-1,1}}\,+\,\sum_{i=1}^{n-1}
  e^{z_{n-1,i}-x_{n,i}}\,+\,g_{i+1}e^{x_{n,i+1}-z_{n-1,i}}\Big)\,\Big\},
 \ee
intertwines the following $D_n$ and $C_{n-1}$ quadratic
Hamiltonians
 \be\label{DC}
  H^{C_{n-1}}(\uz_{n-1})\cdot
  Q^{D_n}_{\,\,\,C_{n-1}}(\ux_n;\,\uz_{n-1})\,=\,
  Q^{D_n}_{\,\,\,C_{n-1}}(\ux_n;\,\uz_{n-1}) H^{D_n}(\ux_n).
 \ee
 To describe a recursion structure of the ground state eigenfunction of
 $C_n$-Toda chain let us introduce the integral operators
$Q^{C_{k}}_{\,\,\,\,C_{k-1}}$ with the kernels given for
$k=2,\ldots,n$  by the convolutions of the intertwining operators
$Q_{C_{k}}^{\,\,\,\,D_{k}}$ and $Q_{D_{k}}^{\,\,\,\,C_{k-1}}$:
 \bqa\label{QC}
  Q^{C_{k}}_{\,\,\,\,C_{k-1}}(\uz_{k};\,\uz_{k-1})\,=\,\int\limits_{\Gamma^{(k)}}
  Q^{C_{k}}_{\,\,\,\,D_{k}}(\uz_{k};\,\ux_{k})\,
  Q^{D_{k}}_{\,\,\,\,C_{k-1}}(\ux_{k};\,\uz_{k-1})\,d\ux_{k}\,.\eqa
  For $k=1$, we set
  \be Q^{C_1}_{\,\,\,\,C_0}(z_{1,1}):=\Psi^{C_1}(z_{1,1})=\int\limits_{\Gamma^{(1)}}
  \exp\{-g_1e^{x_{1,1}+z_{1,1}}-e^{z_{1,1}-x_{1,1}}\}dx_{1,1}.\ee

\npa   For  $k=2,\ldots,n$  introduce the integral operators
$Q^{D_{k}}_{\,\,\,\,D_{k-1}}$ with kernels  given by the
convolutions of the kernels of the intertwining operators
$Q^{D_{k}}_{\,\,\,\,C_{k-1}}$ and $Q^{C_{k-1}}_{\,\,\,\,D_{k-1}}$
 \bqa
  Q^{D_{k}}_{\,\,\,\,D_{k-1}}(\ux_{k};\,\ux_{k-1})\,=\,\int\limits_{\Gamma^{(k)}}
  Q^{D_{k}}_{\,\,\,\,C_{k-1}}(\ux_{k};\,\uz_{k-1})\,
  Q^{C_{k-1}}_{\,\,\,\,D_{k-1}}(\uz_{k-1};\,\ux_{k-1})\,d\uz_{k-1}\,.
 \eqa
 For $k=1$ we set $Q^{D_1}_{\,\,\,\,D_0}(x_{1,1})=1$.

\npa Denote  $\underline{x}_k=(x_{k,1},\ldots,x_{k,k})$,
$\underline{z}_k=(z_{k,1},\ldots,z_{k,k})$. Let us introduce
integral operators $Q^{B_{k}}_{\,\,\,\,B_{k-1}}$ for $k=2,\ldots,n$
with the kernels
 given by the convolutions of the kernels of
elementary intertwiners $Q^{B_{k}}_{\,\,\,\,BC^*_{k}}$ and
$Q^{BC^*_{k}}_{\,\,\,\,B_{k-1}}$:
 \bqa\label{Bkk1}
  Q^{B_{k}}_{\,\,\,\,B_{k-1}}(\ux_{k};\,\ux_{k-1})\,=\,\int\limits_{\Gamma^{(k)}}
  Q^{B_{k}}_{\,\,\,\,BC^*_{k}}(\ux_{k};\, \uz_{k})\,
  Q^{BC^*_{k}}_{\,\,\,\,B_{k-1}} (\uz_{k};\ux_{k-1})\,d\uz_{k}\,.
 \eqa
 For $k=1$ we set
 $$
 Q^{B_1}_{\,\,\,\,B_0}(x_{1,1}):=\Psi^{B_1}(x_{1,1})=
 \int\limits_{\Gamma^{(1)}}\exp\{-2g_1e^{z_{1,1}}-e^{x_{1,1}-z_{1,1}}\}dz_{1,1}.$$

\subsection{Integral representations of Toda chain eigenfunctions}

Using the integral intertwining operators discussed above we can
construct recursively ground state eigenfunctions for quantum
Toda chains associated with classical Lie algebras.
Thus obtained integral representations are compatible with the
representation theory construction of the of generic eigenfunctions
for classical Lie algebras \cite{GLO2}.
The eigenfunction property with respect to quadratic Hamiltonians
follows from the intertwining relations considered in the previous
Subsection. For higher quantum Hamiltonians we should at the moment relay
on the results obtained in representation theory approach.

\begin{prop}
The ground state eigenfunction for $B_n$-Toda chain has the following
integral representation
 \be\label{waveb}
  \Psi^{B_n}(x_1,\ldots,x_n)\,= \int\limits_{\Gamma}\prod_{k=1}^{n}
  Q^{B_{k}}_{\,\,\,\,B_{k-1}}(\ux_{k};
  \,\ux_{k-1})\,\prod\limits_{k=1}^{n-1}d\ux_k,
 \ee
where $x_i:=x_{n,i}$ and $d\ux_k=dx_{k1}\cdots dx_{kk}$.
 The eigenfunction  is given by the integral over a non-compact domain
%$\Gamma=\Gamma_1 \times \prod_{k=1}^{n}\Gamma_2^{(k)}\subset G/B $
$\Gamma$ of the real dimension equal to the complex dimension of the generalized
flag manifold $X=G/B$, where $G=SO(2n+1,\mathbb{C})$ and $B$ is a
Borel subgroup
 \be
\sum_{k=1}^n (2k-1)=n^2=|R_+|. \ee The domain of integration should
be chosen in such a way that the integrand decreases exponentially
at infinities.
\end{prop}

\begin{prop}
The ground state for $C_n$-Toda chain has the following integral
representation
 \be\label{wavecn}
  \Psi^{C_n}(z_1,\ldots,z_n)\,=
  \int\limits_{\Gamma}\prod_{k=1}^{n}
  Q^{C_{k}}_{\,\,\,\,C_{k-1}}(\uz_{k};\,\uz_{k-1})\,\prod\limits_{k=1}^{n-1}d\uz_k\,,
 \ee
where $z_i:=z_{n,i}$ and $d\uz_k=dz_{k1}\cdots dz_{kk}$.
The eigenfunction  is given by the integral
over a non-compact domain
%$\Gamma=\Gamma_1 \times \prod_{k=1}^{n}\Gamma_2^{(k)}\subset G/B$ 
$\Gamma$ of the real dimension
equal to complex dimension of the generalized flag manifold $X=G/B$,
where $G=SO(2n+1,\mathbb{C})$ and $B$ is a Borel subgroup
 \be
\sum_{k=1}^n (2k-1)=n^2=|R_+|. \ee The domain of integration should
be chosen in such a way that the integrand decreases exponentially
at infinities.
\end{prop}

\begin{prop}The ground state eigenfunction for $D_n$-Toda chain
has the following integral representation
 \be\label{WDNone}
  \Psi^{D_n}(x_1,\ldots,x_n)\,=
 \int\limits_{\Gamma}\prod_{k=1}^{n}
  Q^{D_{k}}_{\,\,\,\,D_{k-1}}(\ux_{k};\,\ux_{k-1})\,\prod\limits_{r=1}^{n-1}d\ux_k\,,
 \ee
where $x_i:=x_{n,i}$.
 The eigenfunction  is given by the integral
over a non-compact domain
%$\Gamma=\Gamma_1 \times \prod_{k=1}^{n-1}\Gamma_2^{(k)}\subset G/B$
$\Gamma$ of the real dimension
equal to complex dimension of the generalized flag manifold $X=G/B$,
where $G=SO(2n,\mathbb{C})$ and $B$ is a Borel subgroup
 \be
\sum_{k=1}^{n-1}\, 2k=n(n-1)=|R_+|. \ee
 The domain of integration should
be chosen in such a way that the integrand decreases exponentially
at  infinities.
\end{prop}

\subsection{Quadratic Hamiltonians of $BC_n$- and $I_n$- Toda chains}

Here we write down quadratic Hamiltonians of the generic $BC_n$- and
$I_n$-Toda chains.

\npa Non-reduced root system of type $BC_n$ can be defined as
 \bqa\label{BCroots}
  \alpha_0=2e_1,\hspace{1cm}\alpha_1=e_1,\hspace{1.5cm}
  \alpha_{i+1}=e_{i+1}-e_i,\hspace{2cm} 1\leq i\leq n-1
 \eqa and the
corresponding Dynkin diagram is
 \bqa\label{RBC2}
 \xymatrix{
  \alpha_0\ar@{=>}[dr]\ar@{==>}[dd] &&&&&\\
  & \alpha_2\ar@{-}[r] & \ldots\ar@{-}[r] &
  \alpha_{n-1}\ar@{-}[r] & \alpha_n\\
  \alpha_1\ar@{<=}[ur] &&&&&
 }
 \eqa where the two vertexes connected by a double dash line
correspond to a   reduced $\alpha_1=e_1$ and non-reduced
 $\alpha_0=2e_1=2\alpha_1$ roots.

The corresponding $BC_n$-Toda chain has the following quadratic
Hamiltonian
 \bqa\label{BCHam}
  H^{BC_n}(\ux)\,=\,-\frac{1}{2}\sum_{i=1}^n
  \frac{\partial^2}{\partial x_i^2}\,+\,
  g_1e^{x_1}\,+\,g_2e^{2x_1}\,+\,
  \sum_{i=1}^{n-1}g_{i+2}e^{x_{i+1}-x_i}\,,
 \eqa
where $\ux=(x_1,\ldots,x_n)$. Note that quadratic Hamiltonian
\eqref{BCstar} of the $BC_n^*$-Toda chain is obtained by the following specialization of
 the coupling constants ${g}_1=-g/2$, $g_2=g^2/2$.

\npa  The quadratic Hamiltonian of the  $I_n$-Toda chain \cite{I} is given by
 \be
  H^{I_n}(\ux)\,=\,-\frac{1}{2}\sum_{i=1}^n
  \frac{\pr^2}{\pr x_i^2}\,+\,
  \frac{\widetilde{g}_1}{\bigl(e^{x_1/2}-e^{-x_1/2}\bigr)^2}\,+\,
  \frac{\widetilde{g}_2}{\bigl(e^{x_1}-e^{-x_1}\bigr)^2}\\
  +\,\widetilde{g}_3e^{x_1+x_2}\,+\,
  \sum_{i=1}^{n-1}\widetilde{g}_{i+2}e^{x_{i+1}-x_i}\,,
 \ee
where $\ux=(x_1,\ldots,x_n)$ and $g_i,\,,\,\,i=1,\ldots,n+1$ are
generic coupling constants.

\npa The following specialization of the coupling constants
$\widetilde{g}_i,\,1\leq i\leq n+1$ will be important in the
following
$$
 \widetilde{g}_1=-\frac{g_1}{\sqrt{2g_2}}\,,\hspace{1.5cm}
 \widetilde{g}_2=2\frac{g_1}{\sqrt{2g_2}}\Big(\frac{g_1}{\sqrt{2g_2}}+1\Big)\,,
 \hspace{1.5cm}
 \widetilde{g}_3=g_3\sqrt{\frac{g_2}{2}}\,,
$$
and $\widetilde{g}_{i+2}=g_{i+2}$ for $1<i\leq n-1$.  Denote the
Hamiltonian with the specialized coupling constants by $H^{I^*_n}$:
\be\label{SpecInozemHamiltonian}
H^{I^*_n}(\uz_n)\,=\,-\frac{1}{2}\sum_{i=1}^n
\frac{\pr^2}{\pr z_{n,i}^2}\,\,-\,\,
\frac{\frac{g_1}{\sqrt{2g_2}}}{\bigl(e^{-z_{n,1}/2}-\,e^{z_{n,1}/2}\bigr)^2}\,+\,
\frac{2\frac{g_1}{\sqrt{2g_2}}\Big(\frac{g_1}{\sqrt{2g_2}}+1\Big)}{(e^{-z_{n,1}}-\,
e^{z_{n,1}}\bigr)^2}\\
+\,g_3\sqrt{\frac{g_2}{2}}\,
\Big(e^{z_{n,1}+z_{n,2}}\,+\,e^{z_{n,2}-z_{n,1}}\Big)\,+\,
\sum_{i=2}^{n-1}g_{i+2}e^{z_{n,i+1}-z_{n,i}}\,.
\ee

\subsection{Integral intertwining operators }

Now we conjecture the $BC_n$ and $I_n$ analogs of the intertwining
integral operators. The intertwining relations with
quadratic quantum Hamiltonians are checked straightforwardly.
Hopefully similar intertwining relations hold for higher quantum Hamiltonians.

\begin{lem}
An integral operator  defined by the following kernel:
 \be
  Q_{\quad I^*_n}^{BC_n}(\ux_n,\,\uz_n)\,=\,
  \left(\frac{1-e^{z_{n,1}}}
  {1+e^{z_{n,1}}}\right)^{\frac{g_1}{\sqrt{2g_2}}}\\
  \cdot\exp\Big\{\,-\Big(\sqrt{\frac{g_2}{2}}\,\Big(e^{x_{n,1}+z_{n,1}}\,
  +\,e^{x_{n,1}-z_{n,1}}\Big)\,+\,
  \sum_{i=2}^ng_{i+1}e^{z_{n,i}-x_{n,\,i-1}}\,+\,e^{x_{n,i}-z_{n,i}}\Big)\Big\},
 \ee
intertwines  quadratic $BC_n$ and $I^*_n$
Hamiltonians (\ref{BCHam}) and (\ref{SpecInozemHamiltonian})
 \be\label{BCInozemtsIntertw}
  H^{BC_n}(\ux_n)\cdot Q_{\quad I^*_n}^{BC_n}(\ux_n,\,\uz_n)\,=\,
  Q_{\quad I^*_n}^{BC_n}(\ux_n,\,\uz_n)\cdot H^{I^*_n}(\uz_n).
 \ee
\end{lem}

 \begin{lem} Let $Q_{\quad I^*_{n+1}}^{BC_n}$ be an integral operator
   with the kernel
 \be
  Q_{\quad I^*_{n+1}}^{BC_n}(\ux_n,\,\uz_{n+1})\,=\,
  \left(\frac{1-e^{z_{n+1,1}}}
  {1+e^{z_{n+1,1}}}\right)^{\frac{g_1}{\sqrt{2g_2}}}\\
  \cdot\exp\Big\{\,-\Big(\sqrt{\frac{g_2}{2}}\,\Big(e^{x_{n,1}+z_{n+1,1}}\,
  +\,e^{x_{n,1}-z_{n+1,1}}\Big)\,+\,
  \sum_{i=2}^n\Big(g_{i+1}e^{z_{n+1,i}-x_{n,\,i-1}}+
  e^{x_{n,i}-z_{n+1,i}}\Big)\\+g_{n+2}\,e^{z_{n+1,\,n+1}-x_{nn}}\Big)\Big\}.
 \ee
The  following intertwining relation between  the quadratic
quantum  $BC_n$ and $I^*_{n+1}$  Hamiltonians holds
 \be\label{BCInozemtsRec}
  H^{BC_n}(\ux_n)\cdot Q_{\quad I^*_{n+1}}^{BC_n}(\ux_n,\,\uz_{n+1})\,=\,
  Q_{\quad I^*_{n+1}}^{BC_n}(\ux_n,\,\uz_{n+1})\cdot{H}^{I^*_{n+1}}(\uz_{n+1}).
 \ee
 \end{lem}
Let us denote the kernel of the inverse integral transformation by
$$
 Q^{I^*_{n+1}}_{\quad BC_n}(\uz_{n+1}\,,\ux_n)\,:=\,
 Q_{\quad I^*_{n+1}}^{BC_n}(\ux_n,\,\uz_{n+1})\,.
$$
For any  $k=2,\ldots,n$  define the integral operators
$Q_{BC_{k}}^{\,\,\,\,BC_{k-1}}$
 as
convolutions of the kernels of elementary intertwiners
$Q^{BC_{k}}_{\,\,\,\,I^*_{k}}$ and $Q^{I^*_{k}}_{\,\,\,\,BC_{k-1}}$:
 \bqa
  Q^{BC_{k}}_{\,\,\,\,BC_{k-1}}(\ux_{k};\,\ux_{k-1})\,=\,\int\limits_{\Gamma^{(k)}}
  Q^{BC_{k}}_{\,\,\,\,I^*_{k}}(\ux_{k};\,\uz_{k})\,
  Q^{I^*_{k}}_{\,\,\,\,BC_{k-1}}(\uz_{k};\,\ux_{k-1})\,d\uz_{k}\,.
 \eqa
  For $k=1$ we set
 \be Q^{BC_1}_{\,\,\,\,BC_0}(x_{1,1})= \int\limits_{\Gamma^{(1)}}\Big(\frac{
1-e^{z_{11}}}{1+e^{z_{11}}}\Big)^{\frac{2g_1}{\sqrt{2g_2}}}
\exp\{-\sqrt{g_2/2}(e^{x_{11}+z_{11}}+e^{x_{11}-z_{11}})\}dz_{11}.
 \ee
Let us define the following quadratic $BC$ and $I$ Hamiltonians
\be
  H^{BC_n}(\ux)\,=\,-\frac{1}{2}\sum_{i=1}^n\frac{\pr^2}{\pr x_i^2}\,+\,
  g_1e^{x_1}\,+\,g_2e^{2x_1}+\sum_{k=1}^{n-1}g_{k+2}e^{x_{k+1}-x_k}\,,
 \ee
\be\label{HIn}
  H^{I_n}(\uz)\,=\,-\frac{1}{2}\sum_{i=1}^n\frac{\pr^2}{\pr z_i^2}\,+\,
  \frac{\tilde{g}_1}{\bigl(e^{-z_1/2}-e^{z_1/2}\bigr)^2}+
  \frac{\tilde{g}_2}{\bigl(e^{-z_1}-e^{z_1}\bigr)^2}\\+
  g_3\sqrt{\frac{g_2}{2}}\bigl(e^{z_1+z_2}+e^{z_2-z_1}\bigr)+
  \sum_{k=2}^{n-1}g_{k+2}e^{z_{k+1}-z_k}\,,
 \ee
where
 \be\label{coupling}
  \tilde{g}_1\,=\,-\frac{(2\imath a+1)g_1}{\sqrt{2g_2}}\,,
  \hspace{1.5cm}
  \tilde{g}_2\,=\,2\bigl(\imath a+\frac{g_1}{\sqrt{2g_2}}\bigr)+
  2\Big(\imath a+\frac{g_1}{\sqrt{2g_2}}\Big)^2\,.
 \ee
Note that the coupling constants $\tilde{g}_1$ and $\tilde{g}_2$
are  generic provided the parameter $a$ is generic.
\begin{prop}
The integral operator with the kernel
 \be
  Q^{ I_n}_{\,\,\,\,BC_n}(\uz_n;\,\ux_n|a)\,=\,
  \bigl(1-e^{2z_{n,1}}\bigr)^{-2\imath a}
  \cdot\exp\Big\{\,\imath a\sum_{i=1}^n(z_{n,i}-x_{n,i})\Big\}
  Q^{ I^*_n}_{\,\,\,\,BC_n}(\uz_n;\,\ux_n).
 \ee
intertwines  the quadratic Hamiltonians of $I_n$ and
$BC_n$-type Toda chains  \be\label{QI1}
H^{I_n}(\uz_{n}) Q^{I_n}_{\,\,\,\,BC_n}(\uz_{n};\,\ux_{n}|a)
\,=\,Q^{ I_n}_{\,\,\,\,BC_n}(\uz_{n};\,\ux_{n}|a)\cdot
H^{BC_n}(\ux_{n}).\ee
\end{prop}

\begin{prop}
The integral operator with the kernel
\be
Q^{I_{n+1}}_{\,\,\,\,BC_n}(\uz_{n+1};\,\ux_n)\,=\,
\bigl(1-e^{2z_{n+1,1}}\bigr)^{-2\imath a}
\cdot\exp\Big\{\,\imath
a(\sum_{i=1}^{n+1}z_{n+1,i}-\sum_{i=1}^{n}x_{n,i})\Big\}
\cdot Q^{ I^*_{n+1}}_{\,\,\,\,BC_n}(\uz_n;\,\ux_n),
\ee
intertwines $H^{I_{n+1}}(\uz_{n+1})-\frac{a^2}{2}$  and
$H^{BC_n}(\ux_n)$ Hamiltonians, that is: \be\label{QI2}
\Big(H^{I_{n+1}}(\uz_{n+1})-\frac{a^2}{2}\Big) Q^{
I_{n+1}}_{\,\,\,\,BC_n}(\uz_{n+1};\,\ux_{n}|a) \,=\,Q^{
I_{n+1}}_{\,\,\,\,BC_n}(\uz_{n+1};\,\ux_n|a)\cdot H^{BC_n}(\ux_{n}).
\ee
\end{prop}
Let us denote the kernel of the inverse integral transformation by
$$
Q^{BC_k}_{\quad I_k}(\ux_k\,,\uz_k|a)\,:=\,
Q^{ I_{k}}_{\,\,\,\,BC_k}(\uz_k,\,\ux_{k}|a)\,.
$$
In the following we introduce three types of iterative kernels.  For
any $k=2,\ldots,n$ define the integral operators
$Q_{BC_{k}}^{\,\,\,\,BC_{k-1}}$ depending on auxiliary parameter
 as convolutions of the kernels of elementary intertwiners
$Q^{BC_{k}}_{\,\,\,\,I_{k}}$ and $Q^{I_{k}}_{\,\,\,\,BC_{k-1}}$:
\bqa
Q^{BC_{k}}_{\,\,\,\,BC_{k-1}}(\ux_{k};\,\ux_{k-1}|\lambda_k)\,=\,
\int\limits_{\Gamma^{(k)}}
Q^{BC_{k}}_{\,\,\,\,I_{k}}(\ux_{k};\,\uz_{k}|\lambda_k)\,
Q^{I_{k}}_{\,\,\,\,BC_{k-1}}(\uz_{k};\,\ux_{k-1}|\lambda_k)\,d\uz_{k}\,.
\eqa
  For $k=1$ we set
\be Q^{BC_1}_{\,\,\,\,BC_0}(x_{1,1}|\lambda_1)\,=\,\\
 \int\limits_{\Gamma^{(1)}}\Big(1-e^{2z_{11}}\Big)^{-2\imath\lambda_1}
e^{\imath\lambda_1 z_{11}}\Big(\frac{
1-e^{z_{11}}}{1+e^{z_{11}}}\Big)^{\frac{2g_1}{\sqrt{2g_2}}}
\exp\{-\sqrt{g_2/2}(e^{x_{11}+z_{11}}+
e^{x_{11}-z_{11}})\}dz_{11}.
 \ee
Let us introduce  for $k=2,\ldots,n$  the integral operators
$Q^{I_{k+1}}_{\,\,\,\,I_k}$
\be
 Q^{I_{k}}_{\,\,\,\,I_{k-1}}(\uz_{k};\,\uz_{k-1})\,=\,\int\limits_{\Gamma^{(k)}}
 Q^{I_{k}}_{\,\,\,\,BC_{k-1}}(\uz_{k};\,\ux_{k-1}|a)\,
 Q^{BC_{k-1}}_{\,\,\,\,I_{k-1}}(\ux_{k-1};\,\uz_{k-1}|a)\,d\ux_{k-1}\,.
\ee
For $k=1$ we take
\be
Q^{I_1}_{\,\,\,\,I_0}(z_{11})=(1+e^{z_{11}})^{-\frac{g_1}{\sqrt{2g_2}}-\imath
a}(1-e^{z_{11}})^{\frac{g_1}{\sqrt{2g_2}}-\imath a}e^{\imath
az_{11}}.\ee
Define for $k=2,\ldots,n$  the integral operators
$\Lambda^{I_{k+1}}_{\,\,\,\,I_k}$ :
\be\label{Lambdak}
\Lambda^{I_k}_{\,\,\,\,I_{k-1}}(\uz_k,\uz_{k-1}|\lambda_k)\,=\\\int\limits_{\Gamma^{(1)}}
 \Big(Q^{I_{k}}_{\,\,\,\,BC_{k}}(\uz_{k};\,\ux_{k}|a)
Q^{BC_{k}}_{\,\,\,\,I_{k}}(\ux_{k};\,\underline{u}_{k}|\lambda_k)\Big)
\Big(
Q^{I_{k}}_{\,\,\,\,BC_{k-1}}(\underline{u}_{k};\,\ux_{k-1}|\lambda_k)
Q^{BC_{k-1}}_{\,\,\,\,I_{k-1}}(\ux_{k-1};\,\underline{u}_{k-1}|\lambda_k)\Big)\\
\Big(
Q^{I_{k-1}}_{\,\,\,\,BC_{k-1}}(\underline{u}_{k-1};\,\ux_{k-1}|\lambda_k)
Q^{BC_{k-1}}_{\,\,\,\,I_{k-1}}(\ux_{k-1};\,\underline{z}_{k-1}|a)\Big)
d\ux_{k}\,d\ux_{k-1}\,
d\underline{u_{k}}\,d\underline{u}_{k-1}. \ee For $k=1$ we define
\be
\Lambda^{I_1}_{\,\,\,\,{I_0}}(z_{11}|\lambda_1)\,=\,\int\,dx_{11}\,du_{11}\,
Q^{I_1}_{\,\,\,\,BC_1}(z_{11},x_{11}|a)Q^{BC_1}_{\,\,\,\,I_1}(x_{11},u_{11}|\lambda_1)
Q^{I_1}_{\,\,\,\,BC_0}(u_{11}|\lambda_1), \ee where
$$
Q^{I_1}_{\,\,\,\,BC_0}(u_{11}|\lambda_1)\,=\,
(1+e^{u_{11}})^{-\frac{g_1}{\sqrt{2g_2}}-\imath
 \lambda_1}(1-e^{u_{11}})^{\frac{g_1}{\sqrt{2g_2}}-\imath\lambda_1 }e^{\imath\lambda_1u_{11}}.
$$

\subsection{Givental representations for the eigenfunctions
of $BC_n$ and $I_n$ Toda chains}

In this Section we conjecture  integral representations of
elementary operators which intertwine $I_n$-Toda chain  Hamiltonians
with generic coupling constants and $BC_n$-Toda chain
Hamiltonians with generic coupling constants.
By using intertwining relations we introduce the
conjectural integral formula for the eigenfunction of $I_n$-Toda chain
with generic coupling constants corresponding to zero
eigenvalues.

\begin{conj}
The ground state  eigenfunction for $BC_n$-Toda chain has the
following integral representation
 \be\label{WDN}
  \Psi^{BC_n}(x_1,\ldots,x_n)\,=\,
  \int\limits_{\Gamma}\prod_{k=1}^{n}
  Q^{BC_{k}}_{\,\,\,\,BC_{k-1}}(\ux_{k};\,\ux_{k-1})\,\prod\limits_{k=1}^{n-1}d\ux_k\,,
 \ee
where $x_i:=x_{n,i}$.
The domain  of integration  $\Gamma$  should be chosen in such a way that the integrand
decreases exponentially  at infinities.
\end{conj}

The conjecture is also supported by various degenerations
considering in the next Section. It is also
easily verified for quadratic Hamiltonians using
the intertwining relations (\ref{BCInozemtsIntertw}), (\ref{BCInozemtsRec}).

\begin{conj}
The eigenfunction for $BC_n$-Toda chain corresponding to
generic  eigenvalues \\ parametrized by
$\underline{\lambda}=(\lambda_1,\ldots ,\lambda_n)$ is given by
 \be\label{WDNtwo}
  \Psi^{BC_n}_{\underline{\lambda}}(x_1,\ldots,x_n)\,=\,
  \int\limits_{\Gamma}\prod_{k=1}^{n}
  Q^{BC_{k}}_{\,\,\,\,BC_{k-1}}(\ux_{k};\,\ux_{k-1}|\lambda_k)
\,\prod\limits_{k=1}^{n-1}d\ux_k\,,
 \ee
where $x_i:=x_{n,i}$. The domain  of integration
%$\Gamma=\Gamma_1\times\prod_{k=1}^{n}\Gamma^{(k)}_2$ 
$\Gamma$  should be
chosen in such a way that the integrand
 decreases exponentially at infinities.
\end{conj}
It is easy to check the  conjecture for quadratic Hamiltonians using
the  following iterative procedure: \be
H^{BC_k}(\ux_k)Q^{BC_k}_{BC_{k-1}}(\ux_k;\ux_{k-1}|\lambda_k)\,=\,
Q^{BC_k}_{BC_{k-1}}(\ux_k;\ux_{k-1}|\lambda_k)
\Big(H^{BC_{k-1}}(\ux_{k-1})+\frac{\lambda^2}{2}\Big).
\ee
 At the last step of the iteration we should check the relation:
\be
H^{BC_1}(x_{11})\Psi^{BC_1}_{\lambda_1}(x_{11})=
\frac{\lambda_1^2}{2}\Psi^{BC_1}_{\lambda_1}(x_{11}),
\ee where \be
\Psi_{\lambda_1}^{BC_1}(x_{11})=Q^{BC_1}_{\,\,\,\,BC_0}(x_{11})\,=\,\\
\int_{\Gamma^{(1)}}\Big(1-e^{2z_{11}}\Big)^{-2\imath\lambda_1}\Big(\frac{
1-e^{z_{11}}}{1+e^{z_{11}}}\Big)^{\frac{2g_1}{\sqrt{2g_2}}}
\exp\{-\sqrt{g_2/2}(e^{x_{11}+z_{11}}+ e^{x_{11}-z_{11}})\}dz_{11}.
\ee
The last statement can be checked directly. Notice that the
particular cases $g_1=0$ and $g_2=0$   corresponding to  $B_n$ and
$C_n$-Toda chains respectively, can be transformed (up to certain
multiplier depending on spectrum) to the corresponding eigenfuctions
with non-zero spectrum given in \cite{GLO2}.

%\begin{conj}
%The eigenfunction for $I_n$-Toda chain with generic coupling
%constants given by (\ref{HIn})- (\ref{coupling})  corresponding to
%the eigenvalue equal to $\frac{na^2}{2}$  is given by
% \be\label{waveIn}
%  \Psi^{I_n}(z_1,\ldots,z_n)\,=\,
%  \int\limits_{\Gamma_1}\prod_{k=1}^{n}
%  Q^{I_{k}}_{\,\,\,\,I_{k-1}}(\uz_{k};\,\uz_{k-1})\,\prod\limits_{k=1}^{n-1}d\uz_k\,,
% \ee
%where $z_i:=z_{n,i}$.  The domain of integration
%$\Gamma=\Gamma_1\times\prod_{k=1}^{n}\Gamma^{(k)}_2$  should be
%chosen in such a way that the integrand
% decreases exponentially at possible boundaries and at infinities.
%\end{conj}
%In particular \be
%H^{I_n}(\uz_n)\Psi^{I_n}(z_1,\ldots,z_n)\,=\,
%\frac{na^2}{2}\Psi^{I_n}(z_1,\ldots,z_n).\ee
%It follows from the intertwining relations (\ref{QI1}),(\ref{QI2}).
%and the simply verified fact that \be
%\Big(H^{I_1}(z_{11})-\frac{a^2}{2}\Big)\Psi^{I_1}(z_{11})=0. \ee

\begin{conj}
The eigenfunction for $I_n$-Toda chain with generic coupling
constants given by (\ref{HIn})- (\ref{coupling})  corresponding to
 eigenvalues parametrized by $\underline{\lambda}=(\lambda_1,\ldots ,\lambda_n)$
has the following integral representation
 \be\label{waveIn}
  \Psi^{I_n}_{\lambda_1,\ldots,\lambda_n}(z_1,\ldots,z_n)\,=\,
  \int\limits_{\Gamma}\prod_{k=1}^{n}
  \Lambda^{I_{k}}_{\,\,\,\,I_{k-1}}(\uz_{k};\,\uz_{k-1}|\lambda_k)\,
\prod\limits_{k=1}^{n-1}d\uz_k\,,
 \ee
where $\Lambda^{I_{k}}_{\,\,\,\,I_{k-1}}(\uz_{k};\,\uz_{k-1}|\lambda_k)$
is defined by \eqref{Lambdak} and $z_i:=z_{n,i}$. The domain of integration
%$\Gamma=\Gamma_1\times\prod_{k=1}^{n}\Gamma^{(k)}_3$ 
$\Gamma$  should be chosen in such a way that the integrand
 decreases exponentially  at infinities.
\end{conj}
This conjecture is easily verified for
quadratic Hamiltonian using iteratively the relations
\be
H^{I_k}(\underline{z}_k)\cdot\Lambda^{I_k}_{\,\,\,\,I_{k-1}}(\uz_k;\uz_{k-1})\,=\,
\Lambda^{I_k}_{\,\,\,\,I_{k-1}}(\uz_k;\uz_{k-1})\cdot
\Big(H^{I_{k-1}}(\uz_{k-1})+\frac{\lambda_k^2}{2}\Big),
\ee
and the identity
$$\Big( H^{I_1}(u_{11})-\frac{\lambda^2}{2}\Big)
Q^{I_1}_{\,\,\,\,BC_0}(u_{11}|\lambda_1)\,=\,0.$$

\npa The various limiting case of the Conjectures above are compatible
with the results presented in \cite{GLO2} and recalled in the previous
Section. Let us consider  a specialization  $a=0$
of the integral kernel $Q^{ I_n}_{BC_n}$
 \be
  Q^{I^*_n}_{\,\,\,\,BC_n}(\uz_n;\,\ux_n):=
  \bigl[Q^{ I_n}_{\,\,\,\,BC_n}(\uz_n;\,\ux_n)\bigr]_{a=0}\,=\,
  \Big(
  \frac{1-e^{z_{n,1}}}{1+e^{z_{n,1}}}\Big)^{\frac{g_1}{\sqrt{2g_2}}}\\
  \cdot\exp\Big\{\,-\sqrt{\frac{g_2}{2}}
  \Big(e^{x_{n,1}+z_{n,1}}+e^{x_{n,1}-z_{n,1}}\Big)\,-\,
  \sum_{k=1}^{n-1}g_{k+2}e^{z_{n,k+1}-x_{n,k}}+e^{x_{n,k+1}-z_{n,k+1}}
  \Big\},
 \ee
that intertwines $BC_n$-Toda chain quadratic Hamiltonian:
 \be
  H^{BC_n}(\ux_n)\,=\,
  -\frac{1}{2}\sum_{i=1}^n\frac{\pr^2}{\pr x_{n,i}^2}\,+\,
  g_1e^{x_{n,1}}\,+\,g_2e^{2x_{n,1}}\,+\,
  \sum_{k=1}^{n-1}g_{k+2}e^{x_{n,k+1}-x_{n,k}}\,,\\
 \ee
 with the specialized $I_n$-Toda chain quadratic Hamiltonian:
 \be
  H^{I^*_n}(\uz_n)\,=\,
  -\frac{1}{2}\sum_{i=1}^n\frac{\pr^2}{\pr z_{n,i}^2}\,+\,
  \frac{\tilde{g}_1}{\bigl(e^{-z_{n,1}/2}-e^{z_{n,1}/2}\bigr)^2}+
  \frac{\tilde{g}_2}{\bigl(e^{-z_{n,1}}-e^{z_{n,1}}\bigr)^2}\\+
  g_3\sqrt{\frac{g_2}{2}}
  \bigl(e^{z_{n,1}+z_{n,2}}+e^{z_{n,2}-z_{n,1}}\bigr)+
  \sum_{k=2}^{n-1}g_{k+2}e^{z_{n,k+1}-z_{n,k}}\,,
 \ee
 where
 $$
  \tilde{g}_1\,=\,-\frac{g_1}{\sqrt{2g_2}}\,,
  \hspace{1.5cm}
  \tilde{g}_2\,=\,2\Big(\frac{g_1}{\sqrt{2g_2}}\Big)\Big(
  \frac{g_1}{\sqrt{2g_2}}+1\Big)\,.
$$
The Hamiltonian $H^{I^*_n}(\uz)$ is not generic since its coupling
constants $\tilde{g}_1$ and $\tilde{g}_2$ depend only on the
combination $\frac{g_1}{\sqrt{2g_2}}$.
Shifting the variable $z_1\to z_1+\frac{1}{2}\ln\frac{g_2}{2}$
 one obtains the following expressions for the specialized
$H^{I^*_n}(\uz_n)$ and $Q^{I^*_n}_{\,\,BC_n}(\uz_n;\ux_n)$ when
$a=0$ (tildes over $H^{I^*_n}$ and $Q_{\,\,I^*_n}^{BC_n}$ below
denote the application of the shift of variable $z_1$):
 \be
  \widetilde{H}^{I^*_n}(\uz_n)\,=\,-\frac{1}{2}\sum_{i=1}^n
  \frac{\pr^2}{\pr z_{n,i}^2}\,\,-\,\,\frac{g_1}{2}
  \frac{e^{-z_{n,1}}+\frac{g_2}{2}e^{z_{n,1}}-g_1}
  {\bigl(e^{-z_{n,1}}-\,\frac{g_2}{2}e^{z_{n,1}}\bigr)^2}
  \,+\,\frac{g_2g_3}{2}\,e^{z_{n,1}+z_{n,2}}\,+\,
  \sum_{i=1}^{n-1}g_{i+2}e^{z_{n,i+1}-z_{n,i}}\,,\\\nonumber
  \widetilde{Q}^{I^*_n}_{\,\,\,\,BC_n}(\uz_n,\,\ux_n)\,=\,
  \left(\frac{1-\sqrt{\frac{g_2}{2}}\,e^{z_{n,1}}}
  {1+\sqrt{\frac{g_2}{2}}\,e^{z_{n,1}}}\right)^{\frac{g_1}{\sqrt{2g_2}}}\\
  \cdot\exp\Big\{\,-\Big(\frac{g_2}{2}e^{x_{n,1}+z_{n,1}}+e^{x_{n,1}-z_{n,1}}\,+\,
  \sum_{k=1}^{n-1}g_{k+2}e^{z_{n,k+1}-x_{n,k}}+e^{x_{n,k+1}-z_{n,k+1}}\Big)\Big\}
 \ee
Then one can directly verify the following limiting relations between
various Toda chains:
 \begin{enumerate}
  \item When $g_1\to0$ one has
  \be
   \lim_{g_1\to0} H^{BC_n}(\ux_n)\,=\,H^{C_n}(\ux_n)\,,
   \hspace{1.5cm}
   \lim_{g_1\to0}
   \widetilde{H}^{I^*_n}(\uz_n)\,=\,H^{D_n}(\uz_n)\,,\\
   \lim_{g_1\to0}\widetilde{Q}^{I^*_n}_{\,\,\,\,BC_n}(\uz_n,\,\ux_n)\,=\,
   Q^{D_n}_{\,\,\,\,C_n}(\uz_n,\,\ux_n)\,.
  \ee
 \item When $g_2\to0$ one obtains
 \be
  \lim_{g_2\to0} H^{BC_n}(\ux_n)\,=\,H^{B_n}(\ux_n)\,,
  \hspace{1.5cm}
  \lim_{g_2\to0}
  \widetilde{H}^{I^*_n}(\uz_n)\,=\,H^{BC^*_n}(\uz_n)\,,\\
  \lim_{g_2\to0}\widetilde{Q}^{I^*_n}_{\,\,\,\,BC_n}(\uz_n,\,\ux_n)\,=\,
  Q^{BC^*_n}_{\,\,\,\,B_n}(\uz_n,\,\ux_n)\,.
 \ee
 \end{enumerate}
The last formula readily follows from a simple observation:
$$
 \lim_{\varepsilon\to0}\frac{g}{\sqrt{2\varepsilon}}\,
 \ln\frac{1-\sqrt{\frac{\varepsilon}{2}}\,e^{z}}
 {1+\sqrt{\frac{\varepsilon}{2}}\,e^{z}}\,=\,
 \lim_{\varepsilon\to0}\frac{g}{\sqrt{2\varepsilon}}\Big(-
 2\sqrt{\frac{\varepsilon}{2}}\,e^{z}+o(\varepsilon)\Big)\,=\,
 -g\,e^z
$$

\section{Intertwiners for affine  Toda chains}

In this section we generalize the construction of the elementary
intertwiners to the classical series of affine Lie algebras. For the
necessary facts in the theory of affine Lie algebras see \cite{K}.
Let us first recall the construction of the $Q$-operator for
$A_n^{(1)}$  Toda chain~\cite{PG}. Using the elementary
intertwiners we propose integral operators intertwining affine Toda chains
corresponding to the affine Lie algebras of the same rank.

\subsection{Quadratic Hamiltonians of affine Toda chains}

Let us start with the a list of quadratic Hamiltonians of affine
Toda chains considered below.

\npa Simple roots for the system of type $A_n^{(1)}$ in the standard
basis $\{e_1,\ldots,e_n\}$ read as follows:
$$
 \alpha_0=e_1-e_{n+1}\,,\hspace{2cm}
 \alpha_i=e_{i+1}-e_i\,,\hspace{1cm}
 1\leq i\leq n\,.
$$
The corresponding Dynkin diagram has the following form:
 \bqa
  \xymatrix{
  && \alpha_0\ar@{-}[dll]\ar@{-}[drr] &&\\
  \alpha_1\ar@{-}[r] & \alpha_2\ar@{-}[r] & \ldots &
  \alpha_{n-1}\ar@{-}[l]\ar@{-}[r] & \alpha_n
  }
 \eqa
Quadratic Hamiltonian operator for the $A_n^{(1)}$ closed Toda chain
reads
 \bqa
  \CH^{A_n^{(1)}}(\ux_{n+1})=-\frac{1}{2}\sum_{i=1}^{n+1}
  \frac{\partial^2}{\partial x_{n+1,i}^2}\,+\,
  g_1e^{x_{n+1,1}-x_{n+1,\,n+1}}+
  \sum_{i=1}^ng_{i+1}e^{x_{n+1,\,i+1}-x_{n+1,i}}\,.
 \eqa

\npa Simple roots  of  the twisted affine root system $A^{(2)}_{2n}$
can be expressed  in terms of the standard basis $\{e_i\}$ as
follows
 \be
  \alpha_1=e_1,\hspace{1.5cm} \alpha_{i+1}=e_{i+1}-e_i,\hspace{2cm}
  1\leq i\leq n-1,\\ \nonumber \alpha_{n+1}=-2e_n.
 \ee
The corresponding Dynkin diagram is given  by
 \bqa
  \xymatrix{
  \alpha_1\ar@{<=}[r] & \alpha_2 & \ldots\ar@{-}[l] &
  \alpha_{n-1}\ar@{-}[l]\ar@{<=}[r] & \alpha_n
  }\,\,\,.
 \eqa
Introduce the affine non-reduced  root system $BC'_{n}$ with simple
roots given by
 \be\label{BCDroots}
  \alpha_0=2e_1,\qquad\alpha_1=e_1,\hspace{1.5cm}
  \alpha_{i+1}=e_{i+1}-e_i,\hspace{2cm}1\leq i\leq n-1,\\ \nonumber
  \alpha_{n+1}=-e_n-e_{n-1},
 \ee
and the corresponding Dynkin diagram is as follows
 \bqa
 \xymatrix{
  \alpha_0\ar@{=>}[dr]\ar@{==>}[dd] &&&& \alpha_{n}\\
  & \alpha_2\ar@{-}[r] & \ldots\ar@{-}[r] &
  \alpha_{n-1}\ar@{-}[ur]\ar@{-}[dr] &\\
  \alpha_1\ar@{<=}[ur] &&&& \alpha_{n+1}
 }
 \eqa

\npa Simple roots  of  the twisted affine root system
$A^{(2)}_{2n-1}$ are given by
 \be
  \alpha_1=2e_1,\qquad \alpha_{i+1}=e_{i+1}-e_i,\qquad
  1\leq i\leq n-1\,,\\ \alpha_{n+1}=-e_n-e_{n-1},
 \ee
and corresponding Dynkin diagram is
 \bqa
 \xymatrix{
  &&&& \alpha_{n}\\
  \alpha_1\ar@{=>}[r] & \alpha_2\ar@{-}[r] & \ldots\ar@{-}[r] &
  \alpha_{n-1}\ar@{-}[ur]\ar@{-}[dr] &\\
  &&&& \alpha_{n+1}
 }
 \eqa

\npa Simple roots of  the  affine root system $B^{(1)}_{n}$ are
given by
 \be
  \alpha_1=e_1,\qquad \alpha_{i+1}=e_{i+1}-e_i,\qquad
  1\leq i\leq n-1\,,\\ \alpha_{n+1}=-e_n-e_{n-1},
 \ee
and corresponding Dynkin diagram is as follows
 \bqa
 \xymatrix{
  &&&& \alpha_{n}\\
  \alpha_1\ar@{<=}[r] & \alpha_2\ar@{-}[r] & \ldots\ar@{-}[r] &
  \alpha_{n-1}\ar@{-}[ur]\ar@{-}[dr] &\\
  &&&& \alpha_{n+1}
 }
 \eqa

\npa Introduce the affine non-reduced root system $BC''_{n}$ with
simple roots given  by
 \be
  \alpha_0=2e_1,\quad\alpha_1=e_1,\qquad
  \alpha_{i+1}=e_{i+1}-e_i,\qquad 1\leq i\leq n-1\,,\\
  \alpha_{n+1}=-2e_n,
 \ee
and corresponding Dynkin diagram is given by
$$
 \xymatrix{
  \alpha_0\ar@{=>}[dr]\ar@{==>}[dd] &&&&&\\
  & \alpha_2\ar@{-}[r] & \ldots\ar@{-}[r] &
  \alpha_{n-1}\ar@{-}[r] & \alpha_n\ar@{<=}[r] & \alpha_{n+1}\\
  \alpha_1\ar@{<=}[ur] &&&&&
 }
$$

\npa Simple roots of  the affine root system $C^{(1)}_{n}$ are
 \be
  \alpha_1=2e_1,\qquad \alpha_{i+1}=e_{i+1}-e_i,\qquad
  1\leq i\leq n-1\,,\\ \alpha_{n+1}=-2e_n,
 \ee
and corresponding Dynkin diagram is given  by
$$
 \xymatrix{
  \alpha_1\ar@{=>}[r] & \alpha_2 & \ldots\ar@{-}[l] &
  \alpha_{n-1}\ar@{-}[l]\ar@{-}[r] & \alpha_n\ar@{<=}[r] &
  \alpha_{n+1}
 }\,.
$$

\npa Simple roots of  the affine root system $D^{(1)}_{n}$ are
 \be
  \alpha_1=e_1+e_2,\qquad \alpha_{i+1}=e_{i+1}-e_i,\qquad
  1\leq i\leq n-1\,,\\ \alpha_{n+1}=-e_n-e_{n-1},
 \ee
and corresponding Dynkin diagram is given  by
$$
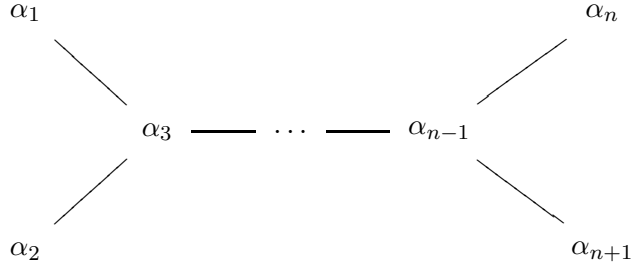
\xymatrix{
 \alpha_1\ar@{-}[dr] &&&& \alpha_n\\
 & \alpha_3\ar@{-}[r] & \ldots\ar@{-}[r] &
 \alpha_{n-1}\ar@{-}[ur]\ar@{-}[dr] &\\
 \alpha_2\ar@{-}[ur] &&&& \alpha_{n+1}
}
$$

\npa Simple roots of  the affine root system $D^{(2)}_n$ are
 \be
  \alpha_1=e_1,\qquad \alpha_i=e_{i+1}-e_i,\qquad 1\leq i\leq
  n-1,\\ \alpha_{n+1}=-e_n,
 \ee and corresponding Dynkin diagram is
given  by
$$
 \xymatrix{
  \alpha_1\ar@{<=}[r] & \alpha_1\ar@{-}[r] & \ldots\ar@{-}[r] &
  \alpha_n\ar@{=>}[r] & \alpha_{n+1}
 }\,.
$$

\npa Simple roots of the (non-reduced) affine root system
$\widehat{BC}_{n+1}$ read as follows
 \be\label{DoubleBCroots}
  \alpha_0=2e_1\qquad\alpha_1=e_1,\\
  \alpha_{i+1}=e_{i+1}-e_i,\hspace{2cm}
  1\leq i\leq n,\\
  \alpha_{n+2}=-e_{n+1}\qquad\alpha_{n+3}=-2e_{n+1},
 \ee and
corresponding Dynkin diagram is given  by
$$
 
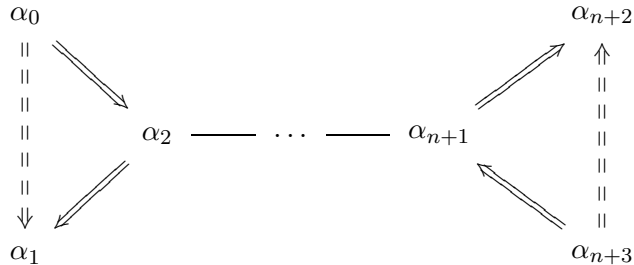
\xymatrix{
  \alpha_0\ar@{=>}[dr]\ar@{==>}[dd] &&&& \alpha_{n+2}\ar@{<=}[dl]\\
  & \alpha_2\ar@{-}[r] & \ldots\ar@{-}[r] &
  \alpha_{n+1}\\
  \alpha_1\ar@{<=}[ur] &&&& \alpha_{n+3}\ar@{=>}[ul]\ar@{==>}[uu]
 }
$$
Let us consider generic quadratic Toda Hamiltonian associated to the root
system $\widehat{BC}_n$ (\ref{DoubleBCroots}):
 \be
  \CH^{\widehat{BC}_n}(\ux_n)\,=\,
 -\frac{1}{2}\sum_{i=1}^n
  \frac{\partial^2}{\partial x_{n,i}^2}\,+\,
  g_1e^{x_{n,1}}+g_2e^{2x_{n,1}}\\
  +\,\sum_{i=1}^{n-1}g_{i+2}e^{x_{n,\,i+1}-x_{n,i}}\,+\,
  g_{n+2}e^{-2x_{nn}}+g_{n+3}e^{-x_{nn}}\,.
 \ee

\subsection{Elementary Intertwiners }

In \cite{PG} Pasquier-Gaudin introduce the  following integral $\CQ$-operator
defined by the integral kernel
 \be
  Q^{A_n^{(1)}}(x_1,\ldots,x_{n+1};\,y_1,\ldots,y_{n+1})\,=\,
  \exp\Big\{\,
  -\sum_{i=1}^{n+1}\Big(e^{x_i-y_i}+g_{i+1}e^{y_{i+1}-x_i}\Big)\,
  \Big\}\,,
 \ee
where $y_{n+2}=y_1, g_{n+2}=g_1$. According to \cite{PG} this integral
kernel satisfies the following intertwining relations
$$
 \CH^{A_n^{(1)}}(\ux_{n+1})\cdot
 Q^{A_n^{(1)}}(\ux_{n+1};\,\uy_{n+1})\,=\,
 Q^{A_n^{(1)}}(\ux_{n+1};\,\uy_{n+1})\cdot
 \CH^{A_n^{(1)}}(\uy_{n+1})\,.
$$
Below we list  kernels of the  integral operators that
conjecturally intertwine
Hamiltonians of affine Toda chains corresponding to  other
classical series of (twisted) affine Lie algebras. We checked
these intertwining relation against the quadratic Hamiltonians.

\npa The integral operator with the following kernel
 \be\label{inttwA2even}
  Q^{A^{(2)}_{2n}}_{\,\,\,\,BC'_{n+1}}(\ux_n,\,\uz_{n+1})\,=\,
  \exp\Big\{-g_1e^{z_{n+1,1}}-\sum_{i=1}^n\Big(
  e^{x_{n,i}-z_{n+1,i}}+g_{i+1}e^{z_{n+1,\,i+1}-x_{n,i}}\Big)\\
  -g_{n+2}e^{-z_{n+1,\,n+1}-x_{nn}}\,\Big\}\,,
 \ee
intertwines Hamiltonian operators for $A_{2n}^{(2)}$- and
$BC_{n+1}^{(2)}$-Toda chains
 \bqa
  \CH^{A_{2n}^{(2)}}(\ux_n)\,=\,-\frac{1}{2}\sum_{i=1}^{n}
  \frac{\partial^2}{\partial x_{n,i}^2}\,+\,g_1e^{x_{n,1}}\,+\,
  \sum_{i=1}^{n-1}g_{i+1}e^{x_{n,i+1}-x_{n,i}}\,+\,
  2g_{n+1}g_{n+2}e^{-2x_{nn}},\\ \nonumber
  \widetilde{\CH}^{BC'_{n+1}}(\uz_{n+1})\,=\,-\frac{1}{2}\sum_{i=1}^{n+1}
  \frac{\partial^2}{\partial z_{n+1,i}^2}\,-\,
  \frac{g_1}{2}e^{z_{n+1,1}}+\frac{g_1^2}{2}e^{2z_{n+1,1}}\\
  \,+\,\sum_{i=1}^{n}g_{i+1}e^{z_{n+1,\,i+1}-z_{n+1,i}}\,+\,
  g_{n+2}e^{-z_{n+1,\,n+1}-z_{n+1,n}}\,.\nonumber
 \eqa
The elementary intertwiner (\ref{inttwA2even}) is involutive; we
use  the following notation for the integral kernel of the inverse
transformation:
 \be
  Q_{\,\,\,\,A^{(2)}_{2n}}^{BC'_{n+1}}(\uz_{n+1},\,\ux_n)\,:=\,
  Q^{A^{(2)}_{2n}}_{\,\,\,\,BC'_{n+1}}(\ux_n,\,\uz_{n+1})\,.
 \ee

\npa The integral operator represented  by the following kernel
 \be\label{inttwA2odd}
  Q^{A^{(2)}_{2n-1}}_{\,\,\,\,A^{(2)}_{2n-1}}(\ux_n,\,\uz_n)\,=\,
  \exp\Big\{-g_1e^{x_{n,1}+z_{n,1}}\\- \sum_{i=1}^{n-1}\Big(
  e^{x_{n,i}-z_{n,i}}+g_{i+1}e^{z_{n,i+1}-x_{n,i}}\Big)-
  e^{x_{nn}-z_{nn}}-g_{n+1}e^{-x_{nn}-z_{nn}}\,\Big\},
 \ee
intertwines Hamiltonian operators for $A^{(2)}_{2n-1}$-Toda
chains with generic coupling constants
 \bqa
  \CH^{A_{2n-1}^{(2)}}(\ux_n)\,=\,-\frac{1}{2}\sum_{i=1}^{n}
  \frac{\partial^2}{\partial x_{n,i}^2}\,+\,
  2g_1e^{2x_{n,1}}+\sum_{i=1}^{n-1}g_{i+1}e^{x_{n,i+1}-x_{n,i}}\\
  \nonumber +\,g_n g_{n+1}e^{-x_{nn}-x_{n,n-1}},\\
  \widetilde{\CH}^{A_{2n-1}^{(2)}}(\uz_n)\,=\,-\frac{1}{2}\sum_{i=1}^{n}
  \frac{\partial^2}{\partial z_{n,i}^2}\,+\,
  g_1g_2e^{z_{n,1}+z_{n,2}}\,+\,
  \sum_{i=1}^{n-1}g_{i+1} e^{z_{n,i+1}-z_{n,i}}\\ \nonumber
  +\,2g_{n+1}e^{-2z_{nn}}\,.
 \eqa
In the following we use a separate notation for integral kernel of the inverse
transformation:
$$
 Q^{A^{(2)}_{2n-1}}_{\,\,\,\,A^{(2)}_{2n-1}}(\uz_n,\,\ux_n)\,:=\,
 Q^{A^{(2)}_{2n-1}}_{\,\,\,\,A^{(2)}_{2n-1}}(\ux_n,\,\uz_n).
$$

\npa  The integral operator represented by the following kernel
 \be\label{inttwBn}
  Q^{B^{(1)}_n}_{\,\,\,\,\,BC''_{n}}(\ux_n,\,\uz_n)\,=\,
  \exp\Big\{-g_1e^{z_{n,1}}\\-\sum_{i=1}^{n-1}\Big(
  e^{x_{n,i}-z_{n,i}}+g_{i+1}e^{z_{n,i+1}-x_{n,i}}\Big)-
  e^{x_{nn}-z_{nn}}-g_{n+1}e^{-x_{nn}-z_{nn}}\,\Big\},
 \ee
intertwines Hamiltonians of $B^{(1)}_n$- and $BC''_{n}$-Toda
chains
 \bqa
  \CH^{B^{(1)}_n}(\ux_n)\,=\,-\frac{1}{2}\sum_{i=1}^{n}
  \frac{\partial^2}{\partial x_{n,i}^2}\,+\,
  g_1e^{x_{n,1}}\,+\,\sum_{i=1}^{n-1}g_{i+1}e^{x_{n,i+1}-x_{n,i}}\\+\,
  g_ng_{n+1}e^{-x_{nn}-x_{n,n-1}},\nonumber \\
  \widetilde{\CH}^{BC''_{n}}(\uz_n)\,=\,-\frac{1}{2}\sum_{i=1}^{n}
  \frac{\partial^2}{\partial z_{n,i}^2}\,-\,
  \frac{g_1}{2}e^{z_{n,1}}+\frac{g_1^2}{2}e^{2z_{n,1}}\\+\,
  \sum_{i=1}^{n-1}g_{i+1}e^{z_{n,i+1}-z_{n,i}}\,+\,2g_{n+1}e^{-2z_{nn}}.\nonumber
 \eqa
 We use the following notation for integral kernel of the inverse
transformation:
 \be
  Q_{\,\,\,\,B^{(1)}_{n}}^{BC''_{n}}(\uz_n,\,\ux_n)\,:=\,
  Q^{B^{(1)}_{n}}_{\,\,\,\,BC''_{n}}(\ux_n,\,\uz_n).
 \ee

\npa The integral operator  with the following kernel
 \be\label{inttwCn}
  Q^{C^{(1)}_{n}}_{\,\,\,\,D^{(1)}_{n+1}}
  (\ux_n,\,\uz_{n+1})\,=\,
  \exp\Big\{-g_1e^{x_{n,1}+z_{n+1,1}}\\-
  \sum_{i=1}^{n}\Big( e^{x_{n,i}-z_{n+1,i}}+
  g_{i+1}e^{z_{n+1,\,i+1}-x_{n,i}}\Big)-
  g_{n+2}e^{-z_{n+1,\,n+1}-x_{nn}}\Big\}
 \ee
intertwines Hamiltonian operators for $C^{(1)}_n$- and
$D^{(1)}_{n+1}$-Toda chains
 \bqa
  \CH^{C^{(1)}_{n}}(\ux_n)\,=\,-\frac{1}{2}\sum_{i=1}^{n}
  \frac{\partial^2}{\partial x_{n,i}^2}\,+\,2g_1e^{2x_{n,1}}\,+\,
  \sum_{i=1}^{n-1}g_{i+1}e^{x_{n,i+1}-x_{n,i}}\,+\,
  2g_{n+1}g_{n+2}e^{-2x_{nn}},\\ \nonumber
  \CH^{D^{(1)}_{n+1}}(\uz_{n+1})\,=\,-\frac{1}{2}\sum_{i=1}^{n+1}
  \frac{\partial^2}{\partial z_{n+1,i}^2}\,+\,
  g_1g_2e^{z_{n+1,1}+z_{n+1,2}}\,+\,
  \sum_{i=1}^{n}g_{i+1}e^{z_{n+1,\,i+1}-z_{n+1,\,i}}\\ +\,
  g_{n+2}e^{-z_{n+1,\,n+1}-z_{n+1,n}}.
 \eqa

\npa The integral operator with the kernel
 \be\label{inttwDn}
  Q^{D^{(1)}_{n}}_{\,\,\,\,C^{(1)}_{n-1}}(\ux_n,\,\uz_{n-1})\,=\,
  \exp\Big\{-g_1e^{x_{n,1}+z_{n-1,1}}\\-\sum_{i=1}^{n-1}\Big(
  e^{z_{n-1,i}-x_{n,i}}+g_{i+1}e^{x_{n,i+1}-z_{n-1,i}}\Big)-
  g_{n+1}e^{-x_{nn}-z_{n-1,\,n-1}}\,\Big\},
 \ee
intertwines Hamiltonian operators for $D^{(1)}_n$-
 and $C^{(1)}_{n-1}$-Toda chains
 \bqa
  \CH^{D^{(1)}_{n}}(\ux_n)\,=\,
  -\frac{1}{2}\sum_{i=1}^{n} \frac{\partial^2}{\partial
  x_{n,i}^2}+g_1g_2e^{x_{n,1}+x_{n,}2}\,+\,
  \sum_{i=1}^{n-1}g_{i+1}e^{x_{n,i+1}-x_{n,i}}\,\\
  \nonumber+\,g_{n+1}e^{-x_{nn}-x_{n,n-1}},\\
  \CH^{C^{(1)}_{n-1}}(\uz_{n-1})\,=\,-\frac{1}{2}\sum_{i=1}^{n-1}
  \frac{\partial^2}{\partial z_{n-1,i}^2}\,+\,
  2g_1e^{2z_{n-1,1}}\,+\,
  \sum_{i=1}^{n-2}g_{i+1}e^{z_{n-1,\,i+1}-z_{n-1,i}}\,\\
  \nonumber+\,2g_{n}g_{n+1}e^{-2z_{n-1,\,n-1}}.
 \eqa
We use the following notations for integral kernels of the inverse
transformations:
 \be
  Q^{D^{(1)}_{n+1}}_{\,\,\,\,C^{(1)}_{n}}
  (\uz_{n+1},\,\ux_n)\,:=\,Q_{\,\,\,\,D^{(1)}_{n+1}}^{C^{(1)}_{n}}
  (\ux_n,\,\uz_{n+1}),\\
  Q_{\,\,\,\,\,D^{(1)}_{n}}^{C^{(1)}_{n-1}}
  (\uz_{n-1},\,\ux_n)\,:=\,Q^{D^{(1)}_{n}}_{\,\,\,\,C^{(1)}_{n-1}}
  (\ux_n,\,\uz_{n-1})\,.
 \ee

\npa The integral operator  with the following kernel
 \be\label{inttwD2n}
  Q^{D^{(2)}_n}_{\,\,\,\,\widehat{BC}^*_{n+1}}
  (\ux_n,\,\uz_{n+1})= \exp\Big\{-g_1e^{z_{n+1,1}}\\-
  \sum_{i=1}^{n}\Big(e^{x_{n,i}-z_{n+1,i}}+
  g_{i+1}e^{z_{n+1,\,i+1}-x_{n,i}}\Big)-
  g_{n+2}e^{-z_{n+1,\,n+1}}\Big\}
 \ee
intertwines Hamiltonian operators for $D^{(2)}_n$- and
$\widehat{BC}_{n+1}$-Toda chains defined as:
 \bqa
  \CH^{D^{(2)}_n}(\ux_n)\,=\,-\frac{1}{2}\sum_{i=1}^n
  \frac{\partial^2}{\partial x_{n,i}^2}\,+\,g_1e^{x_1}\,+\,
  \sum_{i=1}^{n-1}g_{i+1}e^{x_{n,i+1}-x_{n,i}}\,+\,
  g_{n+1}e^{-x_{nn}},\\
  {\CH}^{\widehat{BC}^*_{n+1}}(\uz_{n+1})\,=\,
 -\frac{1}{2}\sum_{i=1}^{n+1}
  \frac{\partial^2}{\partial z_{n+1,i}^2}\,-\,
  \frac{g_1}{2}e^{z_{n+1,1}}+\frac{g_1^2}{2}e^{2z_{n+1,1}}\\
  \nonumber
  +\,\sum_{i=1}^ng_{i+1}e^{z_{n+1,\,i+1}-z_{n+1,i}}\,-\,
  \frac{g_{n+2}}{2}e^{-z_{n+1,\,n+1}}+
  \frac{g_{n+2}^2}{2}e^{-2z_{n+1,\,n+1}}\,.
 \eqa
The integral kernel for the inverse transformation is denoted by
 \be
  Q_{\,\,\,\,D^{(2)}_{n}}^{\widehat{BC}^*_{n+1}}
  (\uz_{n+1},\,\ux_n)\,:=\,
  Q^{D^{(2)}_{n}}_{\,\,\,\,\widehat{BC}^*_{n+1}}
  (\ux_n,\,\uz_{n+1}).
 \ee

\begin{prop} Integral operator with the kernel
 \be
  Q_{\,\,\,\,\widehat{I}_{n+1}}^{\widehat{BC}_n}(\uz_{n+1};\,\ux_n)\,=\,
  \Big(\frac{1-e^{z_{n+1,1}}}{1+e^{z_{n+1,1}}}\Big)^{\frac{g_1}{\sqrt{2g_2}}}
  \Big(\frac{1-e^{z_{n+1,\,n+1}}}{1+e^{z_{n+1,\,n+1}}}
  \Big)^{\frac{g_{n+3}}{\sqrt{2g_{n+2}}}}
  \Big(\frac{1-e^{-2z_{n+1,\,n+1}}}{1-e^{2z_{n+1,1}}}\Big)^a\\
  \cdot\exp\Big\{\,a\Big(\sum_{i=1}^{n+1}z_{n+1,i}\,-\,
  \sum_{k=1}^nx_{n,k}\Big)\,-\,\sqrt{\frac{g_2}{2}}\Big(
  e^{x_{n,1}+z_{n+1,1}}+e^{x_{n,1}-z_{n+1,1}}\Big)\\-\,
  \sum_{k=1}^{n-1}g_{k+2}
  e^{z_{n+1,\,k+1}-x_{n,k}}+e^{x_{n,k+1}-z_{n,k}}\\
  -\sqrt{\frac{g_{n+2}}{2}}\Big(
  e^{z_{n+1,\,n+1}-x_{nn}}+e^{-z_{n+1,\,n+1}-x_{nn}}\Big)
  \,\Big\},
 \ee
intertwines  quadratic Hamiltonians of the $BC_n$-Toda chain:
 \be
  H^{\widehat{BC}_n}(\ux_n)\,=\,
  -\frac{1}{2}\sum_{i=1}^n\frac{\pr^2}{\pr x_{n,i}^2}\,+\,
  g_1e^{x_{n,1}}\,+\,g_2e^{2x_{n,1}}+
  \sum_{k=1}^{n-1}g_{k+2}e^{x_{n,\,k+1}-x_{n,k}}\\
  \nonumber\,+\,
  g_{n+2}e^{-2x_{nn}}-g_{n+3}e^{-x_{nn}}\,,\\
 \ee
 with the closed Inozemtsev quadratic Hamiltonian:
 \be
  H^{\widehat{I}_{n+1}}(\uz_{n+1})\,=\,
  -\frac{1}{2}\sum_{i=1}^n\frac{\pr^2}{\pr z_{n+1,i}^2}
  \,+\,\frac{\tilde{g}_1}
  {\bigl(e^{-z_{n+1,1}/2}-e^{z_{n+1,1}/2}\bigr)^2}\,+\,
  \frac{\tilde{g}_2}{\bigl(e^{-z_{n+1,1}}-e^{z_{n+1,1}}\bigr)^2}\\+
  g_3\sqrt{\frac{g_2}{2}}\bigl(e^{z_{n+1,1}+z_{n+1,2}}+
  e^{z_{n+1,2}-z_{n+1,1}}\bigr)\,+\,
  \sum_{k=2}^{n-1}g_{k+2}e^{z_{n+1,\,k+1}-z_{n+1,k}}\\
  +\sqrt{\frac{g_{n+2}}{2}}\bigl(e^{z_{n+1,1}+z_{n+1,2}}+
  e^{z_{n+1,2}-z_{n+1,1}}\bigr)\,+\,\frac{a^2}{2}\\
  +\frac{\tilde{g}_{n+3}}
  {\bigl(e^{-z_{n+1,\,n+1}}-e^{z_{n+1,\,n+1}}\bigr)^2}\,+\,
  \frac{\tilde{g}_{n+4}}
  {\bigl(e^{-z_{n+1,\,n+1}/2}-e^{z_{n+1,\,n+1}/2}\bigr)^2}\,,
 \ee
where
 \be
  \tilde{g}_1\,=\,-\frac{(2a+1)g_1}{\sqrt{2g_2}}\,,
  \hspace{2.5cm}
  \tilde{g}_2\,=\,2\bigl(a+\frac{g_1}{\sqrt{2g_2}}\bigr)+
  2\Big(a+\frac{g_1}{\sqrt{2g_2}}\Big)^2\,,\\
  \tilde{g}_{n+3}\,=\,
  2\Big(a+\frac{g_{n+3}}{\sqrt{2g_{n+2}}}\Big)^2-
  2\bigl(a-\frac{g_{n+3}}{\sqrt{2g_{n+2}}}\bigr)\,,
  \hspace{1.5cm}
  \tilde{g}_{n+4}\,=\,(2a-1)\frac{g_{n+3}}{\sqrt{2g_{n+2}}}\,.
 \ee
\end{prop}

\subsection{ $Q$-operators for affine Toda chains}

Now we apply the results presented in the previous Sections to  the
construction of the  integral $Q$-operators for all classical
series  of affine Lie algebras as well as affine versions of $BC_n$-
and $I_n$-Toda chains. Let us note that the elementary intertwining
operators and recursive operators for finite Lie algebras can be
obtained from the elementary intertwining operators and
$\CQ$-operators for affine Lie algebras by taking appropriate limits
$g_i\rightarrow0$ in (\ref{inttwA2even})-(\ref{inttwDn}). This
generalizes  known relation between $\CQ$-operator for $A_n^{(1)}$
and recursive operators for $A_n$.
\begin{conj} The integral kernels
for $Q$-operators have the following form
 \bqa
  Q^{A^{(2)}_{2n}}(\ux;\,\uy)\,=\,\int
  Q^{A^{(2)}_{2n}}_{\,\,\,\,BC'_{n+1}}(\ux;\,\uz_{n+1})\,
  Q^{BC'_{n+1}}_{\,\,\,\,A^{(2)}_{2n}}(\uz_{n+1};\,\uy)d\uz_{n+1},
 \eqa
 \bqa
  Q^{A^{(2)}_{2n-1}}(\ux;\,\uy)\,=\,\int
  Q^{A^{(2)}_{2n-1}}_{\,\,\,\,A^{(2)}_{2n-1}}(\ux;\,\uz_n)\,
  Q^{A^{(2)}_{2n-1}}_{\,\,\,\,A^{(2)}_{2n-1}}(\uz_n\,\uy;)\,d\uz_n,
 \eqa
 \bqa
  Q^{B^{(1)}_{n}}(\ux;\,\uy)\,=\,\int
  Q^{B^{(1)}_{n}}_{\,\,\,\,BC''_{n}}(\ux;\,\uz_n)\,
  Q^{BC''_{n}}_{\,\,\,\,B^{(1)}_{n}}(\uz_n;\,\uy)\,d\uz_n\,,
 \eqa
 \bqa
  Q^{C^{(1)}_{n}}(\ux;\,\uy)\,=\,\int
  Q^{C^{(1)}_{n}}_{\,\,\,\,D^{(1)}_{n+1}}(\ux;\,\uz_{n+1})
  Q^{D^{(1)}_{n+1}}_{\,\,\,\,C^{(1)}_{n}}(\uz_{n+1};\,\uy)\,
  d\uz_{n+1}\,,
 \eqa
 \bqa
  Q^{D^{(1)}_{n}}(\ux;\,\uy)\,=\,\int
  Q^{D^{(1)}_{n}}_{\,\,\,\,C^{(1)}_{n-1}}(\ux;\,\uz_{n-1})
  Q^{C^{(1)}_{n-1}}_{\,\,\,\,D^{(1)}_{n}}(\uz_{n-1};\,\uy)\,
  d\uz_{n-1}\,,
 \eqa
 \bqa
  Q^{D^{(2)}_{n}}(\ux;\,\uy)\,=\,\int
  Q^{D^{(2)}_{n}}_{\,\,\,\,\widehat{BC}_{n+1}}(\ux;\,\uz_{n+1})
  Q^{\widehat{BC}_{n+1}}_{\,\,\,\,D^{(2)}_n}(\uz_{n+1};\,\uy)\,
  d\uz_{n+1}\,.
 \eqa
 \bqa
  Q^{\widehat{BC}_n}(\ux;\,\uy)\,=\,\int
  Q^{\widehat{BC}_n}_{\,\,\,\,\widehat{I}_{n+1}}(\ux;\,\uz_{n+1})
  Q^{\widehat{I}_{n+1}}_{\,\,\,\,\widehat{BC}_n}(\uz_{n+1};\,\uy)\,
  d\uz_{n+1}\,,
 \eqa
 \bqa
  Q^{\widehat{I}_n}(\ux;\,\uy)\,=\,\int
  Q_{\,\,\,\,\widehat{BC}_{n-1}}^{\widehat{I}_n}(\ux;\,\uz_{n-1})
  Q_{\,\,\,\,\widehat{I}_n}^{\widehat{BC}_{n-1}}(\uz_{n-1};\,\uy)\,
  d\uz_{n-1}\,,
 \eqa
where $\ux=(x_1,\ldots,x_n)$, $\uy=(y_1,\ldots,y_n)$, and
$\uz_k=(z_{k,1},\ldots,z_{kk})$ for $k=\,n,\,n\pm1$.

\end{conj}

\section{Intertwining operators semi-infinite Toda chains}

Similar approach can be applied to construct $Q$-operators
for infinite root systems $B_{\infty}$, $C_{\infty}$, $BC_{\infty}$
and $D_{\infty}$.

\subsection{Quadratic Hamiltonians for semi-infinite Toda chains}

Simple roots and Dynkin diagrams for infinite Lie algebras
$A_{\infty}$, $B_{\infty}$, $C_{\infty}$, $BC_{\infty}$  and
$D_{\infty}$ are as follows \be
A_{\infty}:\,\,\,\,\,\,\,\,\,\,\,\,\,\qquad\qquad\,\, \qquad
\alpha_{i+1}=e_{i+1}-e_i, \qquad \quad i\in\mathbb{Z}, \ee
$$
 \xymatrix{
  \ldots\ar@{-}[r] & \alpha_{-1}\ar@{-}[r] & \alpha_0 &
  \alpha_{1}\ar@{-}[l]\ar@{-}[r] & \alpha_2\ar@{-}[r] & \ldots
 }\,.
$$

\be B_{\infty}:\,\,\,\,\,\,\,\, \alpha_1=e_1,\qquad
\alpha_{i+1}=e_{i+1}-e_i, \qquad\qquad i\in\mathbb{Z}_{>0}, \ee
$$
 \xymatrix{
  \alpha_1\ar@{<=}[r] & \alpha_2\ar@{-}[r] & \alpha_3\ar@{-}[r] &
  \ldots
 }\,.
$$
\be C_{\infty}:\,\,\,\,\,\,\, \alpha_1=2e_1,\qquad
\alpha_{i+1}=e_{i+1}-e_i, \qquad \qquad i\in\mathbb{Z}_{>0}, \ee
$$
 \xymatrix{
  \alpha_1\ar@{=>}[r] & \alpha_2\ar@{-}[r] & \alpha_3\ar@{-}[r] &
  \ldots
 }\,.
$$
\be D_\infty:\,\,\,\,\,\,\ \alpha_1=e_1+e_2,\qquad
\alpha_{i+1}=e_{i+1}-e_i, \qquad\quad i\in\mathbb{Z}_{>0}, \ee
$$
\xymatrix{
 \alpha_1\ar@{-}[dr] &&&\\
 & \alpha_3\ar@{-}[r] & \alpha_4\ar@{-}[r] &
 \ldots\\
 \alpha_2\ar@{-}[ur] &&&
}\,.
$$

\be BC_{\infty}:\,\,\,\,\,\, \alpha_0=2e_1,\qquad
\alpha_1=e_1,\qquad \alpha_{i+1}=e_{i+1}-e_i, \qquad
i\in\mathbb{Z}_{>0}, \ee
$$
 \xymatrix{
  \alpha_0\ar@{=>}[dr]\ar@{==>}[dd] &&&&&\\
  & \alpha_2\ar@{-}[r] & \alpha_3\ar@{-}[r] &
  \ldots\\
  \alpha_1\ar@{<=}[ur] &&&&&
 }\,.
$$

\subsection{Elementary intertwiners for semi-infinite Toda chains}

$Q$-operator for $A_\infty$ infinite Toda chain is known
(see \cite {T} for the classical limit) and has an elementary
structure. Namely, it is given by an integral operator with the
kernel
 \be
  Q(\ux,\,\uy)=\exp\Big\{\,-\sum_{i\in\mathbb{Z}}
  \Big(e^{x_i-y_i}+g_ie^{y_{i+1}-x_i}\Big)\,\Big\},
 \ee
and intertwines $A_{\infty}$ Toda chain Hamiltonians e.g.
 \bqa
  \CH^{A_{\infty}}(\ux)&=&-\frac{1}{2}\sum_{i\in\mathbb{Z}}
  \frac{\partial^2}{\partial x_i^2}\,+\,
  \sum_{i\in \mathbb{Z}} g_i e^{x_{i+1}-x_i},\\
  \CH^{A_{\infty}}(\uy)&=&-\frac{1}{2}\sum_{i\in\mathbb{Z}}
  \frac{\partial^2}{\partial y_i^2}\,+\,
  \sum_{i\in \mathbb{Z}}g_ie^{y_{i+1}-y_i}.
 \eqa
 In this section  the notations $\ux=(x_1,x_2,\ldots)$,
$\uy=(y_1,y_2,\ldots)$ are used.
Below we conjecture a
 generalization of this integral operator
to  other classical series as well as $BC$ and $I$ series of Toda
chains. The intertwining relations with quadratic Hamiltonians
can be  checked straightforwardly.

\npa The integral operator with the kernel
 \be
  Q^{B_\infty}_{\,\,\,\,\,BC^*_\infty}(\ux,\,\uz)\,=\,\exp\Big\{\,
  -g_1e^{z_1}-\sum_{i>0}\left(e^{x_i-z_i}+g_{i+1}e^{z_{i+1}-x_i}\right)\,
  \Big\},
 \ee
intertwines $B_\infty$ and $BC^*_\infty$ Toda chain Hamiltonian
operators
 \bqa
  \CH^{B_{\infty}}(\ux)&=&-\frac{1}{2}\sum_{i=1}^\infty
  \frac{\partial^2}{\partial x_i^2}
  \,+\,g_1e^{x_1}\,+\,\sum_{i=1}^\infty g_{i+1}e^{x_{i+1}-x_i},\\
  {\CH}^{BC^*_{\infty}}(\uz)&=&-\frac{1}{2}\sum_{i=1}^\infty
  \frac{\partial^2}{\partial z_i^2}\,-\,
  \frac{g_1}{2}e^{z_1}+\frac{g_1^2}{2}e^{2z_1}\,+\,
  \sum_{i=1}^\infty g_{i+1}e^{z_{i+1}-z_i}.
 \eqa

\npa Similarly the integral operator with the kernel
 \be
  Q^{C_\infty}_{\,\,\,\,D_\infty}(\ux,\,\uz)\,=\,
  \exp\Big\{-g_1e^{x_1+z_1}-
  \sum_{i>0}\Big(e^{x_i-z_i}+g_{i+1}e^{z_{i+1}-x_i}\Big)\,\Big\},
 \ee
intertwines $C_\infty$ and $D_\infty$  Toda chain Hamiltonian
operators
 \bqa
  \CH^{C_{\infty}}(\ux)&=&-\frac{1}{2}\sum_{i=1}^\infty
  \frac{\partial^2}{\partial x_i^2}\,+\,
  2g_1e^{2x_1}\,+\,\sum_{i=1}^\infty g_{i+1}e^{x_{i+1}-x_i},\\
  \CH^{D_{\infty}}(\uz)&=&-\frac{1}{2}\sum_{i=1}^\infty
  \frac{\partial^2}{\partial z_i^2}
  \,+\,g_1g_2e^{z_1+z_2}\,+\,
  \sum_{i=1}^\infty g_{i+1}e^{z_{i+1}-z_i}\,.
 \eqa
Thus integral $Q$-operator for  $C_\infty$ has as the  kernel
$$
 Q^{C_\infty}(\ux;\uy)\,=\,\int
 Q^{C_\infty}_{\,\,\,\,\,D_\infty}(\ux,\,\uz) \cdot
 Q^{D_{\infty}}_{\,\,\,\,C_{\infty}}(\uz,\,\uy)\,
 \prod_{i=1}^\infty dz_i\,.
$$

\npa The integral operator with the kernel
 \be
  Q^{D_\infty}_{\,\,\,\,C_\infty}(\uz,\,\ux)\,=\,
  \exp\Big\{-g_1e^{x_1+z_1}-\sum_{i>0}\Big(
  e^{z_i-x_i}+g_{i+1}e^{x_{i+1}-z_i}\Big)\,\Big\},
 \ee
intertwines $D_{\infty}$ and $C_{\infty}$ Toda chain Hamiltonian
operators
 \bqa
  \CH^{D_{\infty}}(\ux)&=&-\frac{1}{2}\sum_{i=1}^\infty
  \frac{\partial^2}{\partial x_i^2}\,+\,
  g_1g_2e^{x_1+x_2}\,+\,
  \sum_{i=1}^\infty g_{i+1}e^{x_{i+1}-x_i},\\
  \CH^{C_{\infty}}(z_i)&=&-\frac{1}{2}\sum_{i=1}^\infty
  \frac{\partial^2}{\partial z_i^2}
  \,+\,2g_1e^{2z_1}\,+\,\sum_{i=1}^\infty g_{i+1}e^{z_{i+1}-z_i}.
 \eqa

\npa Let us  consider the quadratic Hamiltonian of $BC_\infty$-Toda chain
with arbitrary coupling constants:
$$
 \CH^{BC_\infty}(\ux)\,=\,-\frac{1}{2}\sum_{i=1}^\infty
 \frac{\partial^2}{\partial x_i^2}\,+\,g_1e^{x_1}\,+\,
 g_2e^{2x_1}\,+\,\sum_{i=1}^\infty g_{i+2}e^{x_{i+1}-x_i}\,.
$$
Then an integral operator with the kernel
 \be
  Q^{BC_\infty}_{\,\,\,\,I_\infty}(\ux;\uz)\,=\,
  \bigl(1+e^{z_{n,1}}\bigr)^{-\frac{g_1}{\sqrt{2g_2}}-a}
  \bigl(1-e^{z_{n,1}}\bigr)^{\frac{g_1}{\sqrt{2g_2}}-a}
  \cdot\exp\Big\{\,a\sum_{i>0}(z_{n,i}-x_{n,i})\\
  -\,\sqrt{\frac{g_2}{2}}\,\Big(e^{x_1+z_1}\,
  +\,e^{x_1-z_1}\Big)\,-\,\sum_{i>0}\Big(
  g_{i+1}e^{z_{i+1}-x_i}\,+\,e^{x_{i+1}-z_{i+1}}\Big)\Big\}
 \ee
intertwines $\CH^{BC_\infty}(\ux)$ with the specified quadratic
Hamiltonian of $I_\infty$ integrable system:
 \be
  {H}^{I_\infty}(\uz)\,=\,-\frac{1}{2}\sum_{i=1}^\infty
  \frac{\pr^2}{\pr z_i^2}\,\,+\,\,
  \frac{\tilde{g}_1}{\bigl(e^{-z_1/2}-e^{z_1/2}\bigr)^2}+
  \frac{\tilde{g}_2}{\bigl(e^{-z_1}-e^{z_1}\bigr)^2}\\
  +\,\sqrt{\frac{g_2}{2}}\,e^{z_1+z_2}\,+\,
  g_3\sqrt{\frac{g_2}{2}}\,e^{z_2-z_1}\,+\,
  \sum_{i=2}^\infty g_{i+2}e^{z_{i+1}-z_i}\,,
 \ee
where the coupling constants $\tilde{g}_1$ and $\tilde{g}_2$ are
given by (\ref{coupling}):
$$
 \tilde{g}_1\,=\,-\frac{(2a+1)g_1}{\sqrt{2g_2}}\,,
 \hspace{1.5cm}
 \tilde{g}_2\,=\,2\bigl(a+\frac{g_1}{\sqrt{2g_2}}\bigr)+
 2\Big(a+\frac{g_1}{\sqrt{2g_2}}\Big)^2\,.
$$

Let us denote by $Q_{\,\,\,\,BC_\infty}^{I_\infty}(\uz;\ux)$ the
kernel of the inverse integral transformation.

\subsection{ $Q$-operators for semi-infinite Toda chains}

Now one can construct $Q$-operators for the semi-infinite Toda
chains, combining the elementary intertwiners obtained in previous
Subsections.

\begin{conj} The  $Q$-operators for semi-infinite Toda chains
are given by the following integral kernels.
 \be
  Q^{B_\infty}(\ux;\uy)\,=\,\int
  Q^{B_\infty}_{\,\,\,\,BC_\infty}(x_i,\,z_i)
  \cdot Q^{BC_\infty}_{\,\,\,\,\,B_\infty}(z_i,\,y_i)\,
  \prod_{i=1}^\infty dz_i\,;
 \ee
 \be
  Q^{C_\infty}(\ux;\uy)\,=\,\int
  Q^{C_\infty}_{\,\,\,\,\,D_\infty}(\ux,\,\uz) \cdot
  Q^{D_{\infty}}_{\,\,\,\,C_{\infty}}(\uz,\,\uy)\,
  \prod_{i=1}^\infty dz_i\,;
 \ee
 \be
  Q^{D_\infty}(\ux;\uy)\,=\,\int
  Q^{D_\infty}_{\,\,\,\,C_\infty}(x_i,\,z_i) \cdot
  Q^{C_\infty}_{\,\,\,\,D_\infty}(z_i,\,y_i)\,
  \prod_{i=1}^\infty dz_i\,;
 \ee
 \be
  Q^{BC_\infty}(\ux;\uy)\,=\,\int
  Q^{BC_\infty}_{\,\,\,\,I_\infty}(x_i,\,z_i) \cdot
  Q^{I_\infty}_{\,\,\,\,BC_\infty}(z_i,\,y_i)\,
  \prod_{i=1}^\infty dz_i\,;
 \ee
 \be
  Q^{I_\infty}(\uz;\uy)\,=\,\int
  Q^{I_\infty}_{\,\,\,\,BC_\infty}(z_i,\,x_i) \cdot
  Q_{\,\,\,\,I_\infty}^{BC_\infty}(x_i,\,y_i)\,
  \prod_{i=1}^\infty dx_i\,.
 \ee
\end{conj}
One can directly check this conjecture for the
intertwining relations with quadratic Hamiltonians.

\vskip 1cm

\noindent {\small {\bf A.G.} {\sl Institute for Theoretical and
Experimental Physics, 117259, Moscow,  Russia; \hspace{8 cm}\,
\hphantom{xxx}  \hspace{2 mm} School of Mathematics, Trinity
College, Dublin 2, Ireland; \hspace{8 cm}\,
\hphantom{xxx}   \hspace{1 mm} Hamilton
Mathematics Institute, Trinity College, Dublin 2, Ireland;}}

\noindent{\small {\bf D.L.} {\sl
 Institute for Theoretical and Experimental Physics,
117259, Moscow, Russia};\\
\hphantom{xxxx} {\it E-mail address}: {\tt lebedev@itep.ru}}\\

\noindent{\small {\bf S.O.} {\sl
 Institute for Theoretical and Experimental Physics,
117259, Moscow, Russia};\\
\hphantom{xxxx} {\it E-mail address}: {\tt Sergey.Oblezin@itep.ru}}

\end{document}